\documentclass[12pt]{article}
\usepackage{amsmath,amssymb}   
\usepackage[mathscr]{euscript}
\usepackage[dvips]{graphicx}
\usepackage{mediabb}
\usepackage{wrapfig}
\usepackage{hangcaption}
\usepackage{color}
\usepackage{amssymb,amsfonts,amsthm}  
\usepackage{bm}  
\usepackage{enumerate}  
\usepackage[dvips]{graphicx}  
\makeatletter

\theoremstyle{plain} 
\newtheorem{thm}{Theorem}[]

\newtheorem{prop}[thm]{Proposition}
\newtheorem{lem}[thm]{Lemma}
\newtheorem{remark}[thm]{Remark}
\newtheorem{definition}[thm]{Definition}
\newenvironment{prf}
   {{\noindent \bf Proof. }}{\hfill \qed}

\newcommand\F{\mathfrak{F}}
\newcommand\R{\mathbb{R}}
\newcommand\cR{\cal{R}}

\newcommand\cQ{\cal{Q}}

\makeatletter
    
    \@addtoreset{equation}{section}
  \makeatother

\begin{document}
\title
 {Pattern formation in a diffusion-ODE model with hysteresis }
\author{Anna Marciniak-Czochra$^1$, Madoka Nakayama$^2$, Izumi Takagi$^3$\\
\vspace{0.3cm}\\
\small $^1$ Institute of Applied Mathematics and BIOQUANT, University of  Heidelberg\\
\small Im Neuenheimer Feld 294, 69120 Heidelberg, GERMANY\\
\small $^2$ Sendai National College of Technology, Natori\\
\small 48 Nodayama, Medeshima-Shiote, Natori, 981-1239, JAPAN \\
\small $^3$ Mathematical Institute, Tohoku  University\\
\small Sendai, 980-8578, JAPAN }

\maketitle

\abstract
Coupling diffusion process of signaling molecules with nonlinear interactions of intracellular processes and cellular growth/transformation 
 leads to  a system of reaction-diffusion equations coupled with ordinary differential equations (diffusion-ODE models), which differ from the usual reaction-diffusion systems. One of the mechanisms of pattern formation in such systems is based on the existence of multiple steady states and hysteresis in the ODE subsystem. Diffusion tries to average different states and is the cause of spatio-temporal patterns. In this paper we provide a systematic description of stationary solutions of such systems, having the form of transition or boundary layers. The solutions are discontinuous in the case of non-diffusing variables whose quasi-stationary dynamics exhibit hysteresis. The considered model is motivated by biological applications and elucidates a possible mechanism of formation of patterns with sharp transitions.

%

\section{Introduction}

One of the most frequently discussed organisms in theoretical papers on biological pattern formation is a fresh-water polyp {\it Hydra}.
{\it Hydra} is a small coelenterate living in fresh water and it is best known for its ability of regeneration.
When its head is cut, 
in a few days a new head completely regenerates.
 The question of de novo pattern formation in a homogenous {\it Hydra} tissue was addressed by Turing in his pioneering paper
\cite{Turing}. Based on Turing's idea,  Gierer and
Meinhardt \cite{G-M} proposed a reaction-diffusion model consisting of an activator and
an inhibitor to explain the regeneration experiment of {\it Hydra}.  Several models have been proposed which modify or refine the activator-inhibitor system by Gierer and Meinhardt, for example by MacWilliams \cite{MacWilliams1982} and Meinhardt \cite{Meinhardt}. 

Another class of mathematical models for pattern formation follows the hypothesis that positional value of the cell is determined by the density of cell-surface receptors, which regulate the expression of genes responsible for cell differentiation. Such models, called receptor-based models,  involve diffusive species and non-diffusive species. The first receptor-based model for {\it Hydra} was  proposed by Sherrat, Maini, J{\"a}ger an M{\"u}ller in \cite{Sherrat}. Later,  receptor-based models without imposing initial gradients were proposed by Marciniak-Czochra \cite{MC03,MC06}.  In general, equations of such models can be represented by the following initial-boundary value  problem:
\begin{align} \label{recmodel}
\left\{
\begin{aligned}
& \frac{\partial U}{\partial t} = D 
 \Delta U +F(U,V) & &\text{for}\quad x\in \Omega, t>0,\\
& \frac{\partial V}{\partial t} = G(U,V) & &
\text{for}\quad  x\in \Omega, t>0,\\
& \frac{\partial U}{\partial \nu} =0  & &
  \text{for}\quad  x\in \partial \Omega,  t>0, \\
& U(x,0) =U_{0}(x), \,\,
V(x,0) =V_{0}(x) & & \text{for} \quad x\in \Omega, 
\end{aligned}
\right.
\end{align} 
where $U$ is a vector of variables describing the dynamics of
diffusing extracellular molecules and enzymes, which provide
cell-to-cell communication, while $V$ is a vector of variables
localized on cells, describing cell surface receptors and
intracellular signaling molecules, transcription factors, mRNA, etc. $F$ and $G$ are  smooth mappings.
$D$ is a diagonal matrix with positive coefficients on the diagonal, $\Omega$ is a bounded domain in $\R^{n}$  with smooth boundary $\partial\Omega$ and $\nu$ denotes the unit outer normal to $\partial\Omega$ . 

A rigorous derivation, using methods of asymptotic analysis
(homogenization) of the macroscopic diffusion-ODE models
describing the interplay between the nonhomogeneous cellular dynamics
and the signaling molecules diffusing in the intercellular space has
been presented in \cite{MCP,MC12}. 

In the framework of reaction-diffusion systems there are essentially two mechanisms of formation of stable spatially heterogenous patterns:
\begin{itemize}
\item diffusion-driven instability (DDI) which leads to destabilization of a spatially homogeneous attractor and emergence of stable spatially heterogenous and spatially regular structures (Turing patterns),
\item a mechanism based on the multistability in the structure of nonnegative spatially homogenous stationary solutions, which leads to transition layer patterns.
\end{itemize} 
DDI and multistability can also coexist yielding different dynamics for different parameter regimes.

The two mechanisms of pattern formation lead to interesting effects in the case of receptor-based models consisting of single reaction-diffusion equation coupled to ordinary differential equations. As shown in \cite{MKS13} for a particular example of a reaction-diffusion-ODE system, it may happen that the system exhibits the DDI but there exist no stable Turing-type patterns  and the emerging spatially heterogenous structures have a dynamical character. In numerical simulations, solutions having the form of periodic or irregular spikes have been observed.
The result on instability of all Turing patterns  can be extended to  general diffusion-ode systems with a single diffusion operator exhibiting DDI.  Consequently, multistability is necessary in such systems to provide stable spatially heterogeneous stationary patterns.
On the other hand, it has been recently shown that the system with multistability but reversible quasi-steady state in the ODE subsystem, i.e. $G(u,v)$ globally invertible, cannot exhibit stable spatially heterogeneous patterns \cite{KMC12}.  Hysteresis is necessary to obtain stable patterns in the diffusion-ODE models with a single diffusion. 

Therefore, in the current work we focus on specific nonlinearities $F$ and $G$ which describe a generalized version of the nonlinearities proposed in \cite{MC06} to model pattern formation in Hydra that exhibit hysteresis-effect. In this case the steady state equation $G(U, V) = 0$ has multiple solutions. The patterns observed in such models are not Turing patterns. In fact, the system does not need to exhibit DDI.  Indeed,  in most cases  its constant steady states do not change stability and spatially heterogenous stationary solutions appear far from equilibrium due to the existence of multiple quasi-steady states.  

Numerical simulations suggest that solutions of the model with the hysteresis-type nonlinearity  behave differently from that of standard reaction-diffusion systems.
For example, some numerical solutions seem to approach quickly steady-states with jump discontinuity.
For a correct understanding of what is actually happening,  it is important to build a rigorous theory on the basic properties of diffusion-ODE systems.

Another important aspect of the hysteresis-based mechanism of pattern formation is related to the co-existence of different steady states.
In particular, bistability in the dynamics of the growth factor controlling cell differentiation in the receptor-based models explains the experimental observations on the multiple head formation in \textit{Hydra}, which is not possible to describe by using Turing-type models \cite{MC06,KMC12}. Those observations showing importance and biological relevance of the rich structure of patterns in diffusion-ODE models have motivated the present work.

In the current paper we study rigorously a certain class of diffusion-ODE systems with hysteresis and derive some of the fundamental properties of solutions such as the boundedness of solutions of the initial-boundary value problem and the existence of initial functions that result in trivial steady-states. 
The novelty of the paper is in providing a systematic description of the stationary solutions of a receptor-based model with hysteresis. The sationary problem corresponding to (1.1) can be reduced to a boundary value problem for a single reaction-diffusion equation with discontinuous nonlinearity. Construction of transition layer solutions for such systems was undertaken in Mimura, Tabata and Hosono \cite{MTH} by using a shooting method. They introduced a diffusion-ODE system as an auxiliary system needed to obtain a steady-state solution with an interior transition layer. The result was applied by Mimura \cite{Mimura} to show the existence of discontinuous patterns in a model with density dependent diffusion. 
While in their models, the transition layer solution was unique, we face the problem of the existence of infinite number of solutions with changing connecting point. To deal with this difficulty we propose a new approach to construct all monotone stationary solutions having either a transition layer or a boundary layer. 

Our results show that the emerging patterns may exhibit discontinuities that may explain sharp transitions in gene expression observed in many biological processes.   In the context of  {\it Hydra} pattern formation, the hysteresis-driven mechanism allows for formation of gradient-like patterns in the expression of Wnt, corresponding to the normal development as well as emergence of patterns with multiple maxima describing transplantation experiments.

The remainder of this paper is organized as follows. Section \ref{Results} presents main results of this paper. 
Section \ref{SectionExistence} is devoted to the existence of mild and H\"older continuous solutions of the initial-boundary value problem \eqref{model},
and  their uniform boundedness. In Section \ref{SectionKinetics} we consider an initial-boundary value problem for the two-component system being a quasi-steady state reduction of the original problem  \eqref{model}. We characterize the regimes of spatially homogeneous dynamics of this model.
Finally, in Section \ref{SectionPattern} we focus on the stationary discontinuous two-point boundary value problem and provide a characterization of its solutions.
We construct a monotone increasing stationary transition layer solution and show how the position of the layer (boundary layer vs interior layer) depends on the model parameters. We discuss also the existence of nonmonotone stationary solutions.


\section{Main results}\label{Results}

\subsection{Statement of the problem}

We consider the following system of equations, defined on a bounded domain $\Omega \subset \R^n$ with a sufficiently smooth boundary $\partial \Omega$:
\begin{align} \label{model}
\left\{
\begin{aligned}
\dfrac{\partial u_{1}}{\partial t} &= -\mu_{1}u_{1}-B(u_1,u_2,u_3)+m_{1}, &&\quad \text{for}\,\, x\in \Omega,\,t>0, \\
 \dfrac{\partial u_{2}}{\partial t} &= -\mu_{2}u_{2}+B(u_1,u_2,u_3),&&\quad \text{for}\,\, x\in \Omega,\,t>0, \\ 
 \dfrac{\partial u_{3}}{\partial t} &=\dfrac{1}{\gamma}\Delta u_{3}-\mu_{3}
 u_{3}-B(u_1,u_2,u_3)+ u_{4},&&\quad \text{for}\,\, x\in \Omega,\,t>0,  \\ 
 \dfrac{\partial u_{4}}{\partial t} &=-\delta u_{4} +P(u_2,u_3,u_4),&& \quad \text{for}\,\, x\in \Omega,\,t>0,\\
  \dfrac{\partial u_{3}}{\partial \nu}&=0&&\quad \text{for}\,\, x \in \partial\Omega,\,t>0,\\
u_{j}(x&,0)=u_{j}^0(x)\quad (j=1,\dots,4)&&\quad \text{for}\,\, x\in \Omega,
\end{aligned}
\right.
\end{align}
where $u_{1}$, $u_{2}$, $u_{3}$ and $u_4$ denote the densities of free
receptors, bound receptors, ligand and transcript of the ligand, respectively. The parameters
 $\mu_{1},\mu_{2},\mu_{3},m_{1},\delta,\gamma$
are positive constants. The initial functions $u^0_{j}(x)$, for $j=1,\dots,4$, are assumed to be smooth and positive on $\overline{\Omega}$.

The function $B(u_1,u_2,u_3)$ describes the process of binding and dissociation.  The simplest example of this function, $B(u_1,u_2,u_3)= b u_1 u_3 - d u_2$ with  $b,d>0$, was considered in  \cite{MC03,MC06, MCP}. The natural decay of all model ingredients as well as the translation process are assumed to be linear. The production of free receptors is constant and the production of ligand transcript is given by a function $P(u_2,u_3, u_4)$. The function $P(u_2,u_3,u_4)$  is modeling the outcome of intracellular signal transduction. In the model proposed in \cite{MC06} it was assumed that the process involves hysteresis-like relation between the signal given by the density of diffusing signaling molecules and the cell response. 

\begin{figure}[t]
 \begin{center}
 \includegraphics[width=60mm]{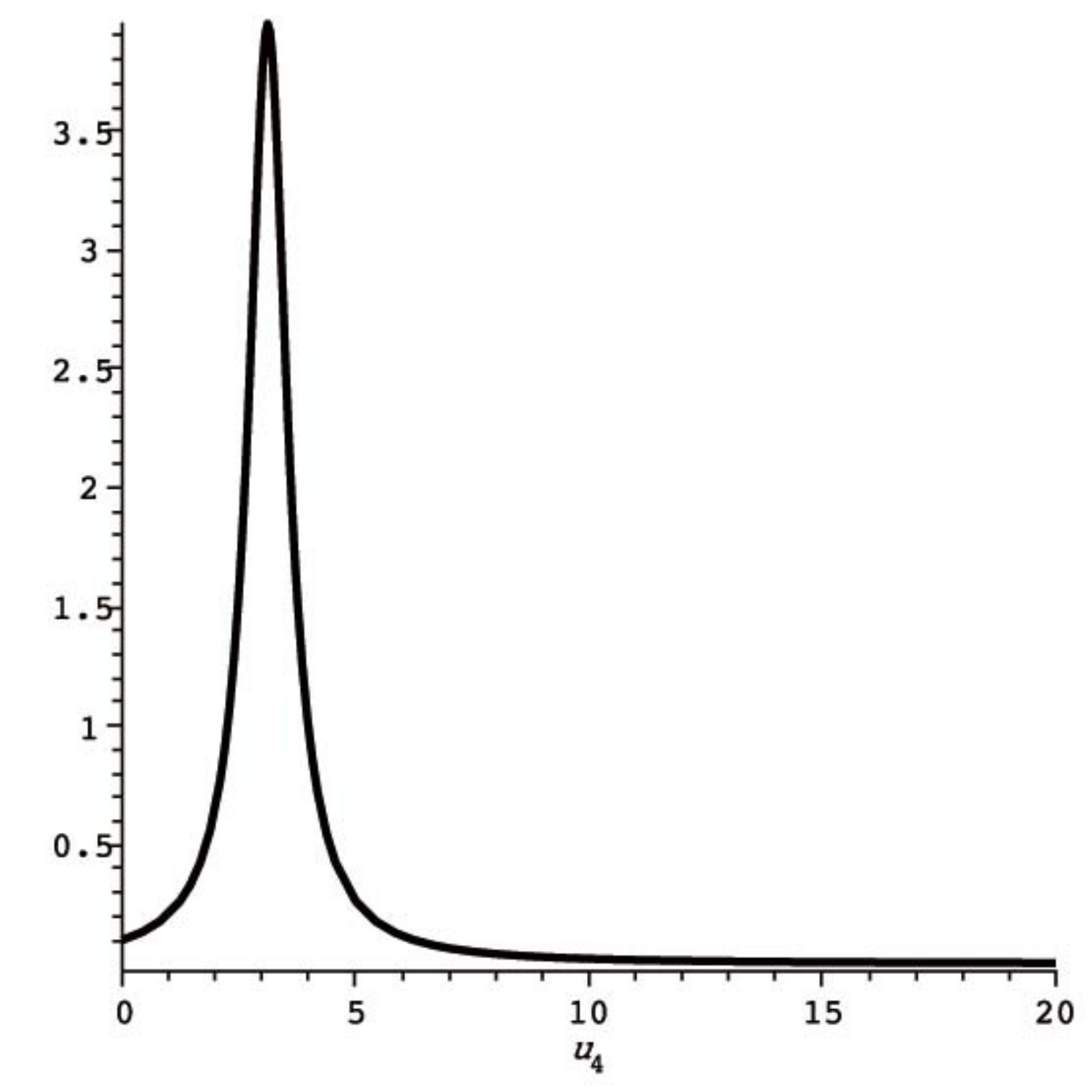} 
 \end{center}
\caption{$S(u_{4})$ when $\sigma=0.1,\,\beta_{l}=0.6244\, \,\text{and}\, \, m_{4}=0.1$  }
 \label{picS}
\end{figure}

In this paper, we  keep the following assumptions concerning the nonlinearities of model \eqref{model}:
\begin{align}
& B(u_1,u_2,u_3) = b u_1 r(u_3) - du_2, \label{AssB}\\
& P(u_2,u_3,u_4) = u_3 S(u_4), \label{AssP}
\end{align}
where $r(u_3)$ is a smooth function modeling binding of ligands to free receptors and $S(u_4)$ is a function modeling a control of the synthesis of new transcripts. We assume $S$ to be a positive function with a maximum at $u_4=u_4^*>0$;  growing for $u_4\in (0,u_4^*)$ and decaying asymptotically to zero for $u_4\rightarrow \infty$, see Fig \ref{picS}.

For computational purposes we take a simple realization of such function in the form 
\begin{equation}\label{AssS}
S= \dfrac{m_4 }{1+\sigma u_{4}^{2}-\beta_{l} u_{4}},
\end{equation}
where  $m_{4},\sigma,\beta_{l},\delta$ are positive constants.
 Throughout the paper we impose the following condition on $\beta_{l}$ and $\sigma$
\begin{align}\label{A1}
\beta_{l}^{2}<4\sigma.
\end{align}
This is a necessary and sufficient condition to provide positivity of the control function $S$ by imposing
\begin{align}\label{1a}
 c_{0}=\min_{u\in\R}(1+\sigma u^{2}-\beta_{l} u)>0.
 \end{align}
 Furthermore, to streamline the presented analysis we chose
 
 \begin{equation}\label{eq_r}
 r(u_3)=u_3.
  \end{equation}
  
\begin{remark}
 {\rm  The specific choice of  nonlinearity $P$ is assumed to make a computation of various quantities easier, yet the new model retains the essence of the mechanism to generate pattern in the original one proposed in \cite{MC06}. In particularly, the direct dependence of $P$ on $u_3$ omitting the signal transduction through $u_2$ is assumed for a mathematical simplicity.  It can be justified by using a quasi-steady state approximation of the receptor dynamics,  assuming additionally that a natural decay of bound receptors is negligible, i.e. $\mu_2=0$. In the case of a positive $\mu_2$ we obtain a more complicated relation, where linear dependence on $u_3$ in \eqref{AssP} should be replaced by a Hill-type function, see the quasi-stationary approximation given in \eqref{QuasiStat}. Such nonlinearity does not change, however, the hysteresis-related properties of the considered system.}
\end{remark}
 

\subsection{Existence and boundedness result}

 First we establish the existence and boundedness of solutions of the initial-boundary value problem \eqref{model}.
 

\begin{thm}\label{ExistBounded}
Let $ \Omega $ be a bounded domain in $\R^n$ with sufficiently smooth boundary $\partial \Omega$. 
Let $u^0_{j}(x)$,  for $j=1,\dots, 4$, be positive and H\"older continuous functions on $\overline \Omega$. 
Suppose moreover that  $u_3^0 \in C^{2+\alpha}(\overline{\Omega}) $, $0<\alpha<1$, and $ \partial u_{3}^{0} / \partial \nu =0$ on $\partial\Omega$.
Then, the initial-boundary value problem \eqref{model} has a unique classical solution $(u_{1}(x,t),\,u_{2}(x,t),\,u_{3}(x,t),$ $u_{4}(x,t))$
for all $t>0$.
Moreover,
there exist positive constants $\rho_{j}\,(j=1,\dots,\,4)$,
depending on the initial functions $(u^0_{1},\,u^0_{2},\,u^0_{3},\,u^0_{4})$,
such that 
$$0<u_{j}(x,t)\leq \rho_{j}\quad \text{for all}\,\, x\in \overline \Omega,\, t\geq0,$$
for $j=1,\,2,\,3,\,4$.
\end{thm}
\begin{remark}
 {\rm 
Some numerical solutions seem to develop a singularity in finite time;
to be more specific,
in some of numerical solutions of \eqref{model},
an interior transition layer is formed and  the spatial derivative of $u_4(x,t)$ becomes larger and larger in the layer,
and $u_4(x,t)$ seems to form a jump discontinuity in finite time.
However, using the regularity of solutions, we can rule out this possibility, as long as we choose the initial data to be sufficiently smooth.}
\end{remark}

\subsection{Reduction of the model}
To reduce the model, we consider the case where the free and bound receptors
are in a quasi-stationary state,
by which we mean that the derivatives in the first two equations are set to zero. Then we obtain
\begin{align}\label{QuasiStat}
u_{1}=\dfrac{m_{1}(d+\mu_{2})}{\mu_{1}(d+\mu_{2})+\mu_{2}b u_{3}}\quad \text{and}\quad u_{2}=\dfrac{m_{1}bu_{3}}{\mu_{1}(d+\mu_{2})+\mu_{2}bu_{3}}.
\end{align}
Substituting \eqref{QuasiStat} in  \eqref{model},  we obtain the following initial-boundary value problem for $u(x,t)=u_{3}(x,t)$ and $v(x,t)=u_{4}(x,t)$:
\begin{align}\label{originalproblem}
\left\{
\begin{aligned}
&\dfrac{\partial u}{\partial t}=\dfrac{1}{\gamma}\Delta u+f(u,v)&&\quad \text{for}\,\,\in \Omega,\,t>0,\\
&\dfrac{\partial v}{\partial t}=g(u,v)&&\quad \text{for}\,\,\in \Omega,\,t>0,\\
&\dfrac{\partial u}{\partial \nu}=0&&\quad \text{for}\,\,x\in \partial\Omega, \,\,t>0,\\
&u(x,0)=u^0_{3}(x),\,v(x,0)=u^0_{4}(x)&&\quad \text{for} \,\,x \in \Omega
\end{aligned}
\right.
\end{align}
with
\begin{align}\label{2eqnonlinearities}
f(u,v)&= v-\mu_{3}u-\dfrac{m_{1}\mu_{2}bu}{\mu_{1}(\mu_{2}+d)+\mu_{2}bu},\\
g(u,v)&= -\delta v+u S(v).
\end{align}
In the remainder of this paper we focus on the reduced model \eqref{originalproblem} defined in one-dimensional domain $\Omega=(0,1)$, except in subsection 2.4.1.


\subsection{Stationary solutions}


\subsubsection{Spatially homogeneous stationary solutions}

We are interested 
in the behavior of the nullclines of the kinetic system 
\begin{align}\label{kineticsystem}
\left\{
\begin{aligned}
&u_t = f(u,v),\\
&v_t= g(u,v).
\end{aligned}
\right.
\end{align}

\begin{prop}
Let \begin{align}\label{A2}
3\sigma<\beta_{l}^{2}<4\sigma.
\end{align}
Then, there exists a range of  $\mu_{3}$ such that
system \eqref{originalproblem}  has three  nonnegative spatially homogeneous stationary solutions: $(u_0,v_0)=(0,0)$, $(u_m,v_m)$ and $(u_1,v_1)$. 
\end{prop}
%
%
\begin{prf}
Obviously,  the trivial solution $(u_{0},v_{0})=(0,0)$ is a stationary solution. 
The nullcline $f(u,v)=0$ defines a strictly monotone increasing function of $u$
$$v=\Phi (u)=\frac{\mu_{3}u+m_{1}\mu_{2}bu}{\mu_{1}(\mu_{2}+d)+\mu_{2}bu}$$ and the nullcline $g(u,v)=0$  defines a cubic function of $v$  $$u=\Psi(v)=\frac{\delta v}{ S(v)}= \frac{\delta v (1+\sigma v^{2}-\beta_{l} v)}{m_{4}},$$ which achieves one positive local maximum at $v=v_{-}$ and one positive local minimum at $v=v_{+}$ if and only if \eqref{A2} is satisfied.

We put $u_{+}=\Psi(v_{-}),\,u_{-}=\Psi(v_{+})$.
Under this condition,
we require further that the two nullclines intersect at exactly two points    
$(u_{m},\,v_{m})$ and $(u_{1},\,v_{1})$ in the first quadrant.
This is made possible by choosing the coefficient $\mu_{3}$ in an appropriate interval:
\begin{align}\label{A3}
\mu_{*}<\mu_{3}<\mu^{*}
\end{align}
where $\mu_{*}$ and $\mu^{*}$ are suitable constants such that $0\leq \mu_{*}<\mu^{*}$.
\end{prf}

%
\begin{lem}\label{saddle}
The equilibrium  $(u_{m},\,v_{m})$ of the kinetic system \eqref{kineticsystem} is a saddle point.
Its stable manifold $W^{s}$ intersects with the positive $u$-axis at $(U_{s},0)$ and with the positive $v$-axis at $(0,V_{s})$.
Put ${\cal{Q}}=\{(u,v)\mid u>0,\,v>0\}$.
The projection of $W^{s}\cap {\cal{Q}}$ onto the $u$-axis coincides with the interval $0<u<U_{s}$ and the projection of $W^{s}\cap {\cal{Q}}$
onto the $v$-axis coincides with the interval $0<v<V_{s}$.
Moreover,
these projections are injective.
\end{lem}

%
\begin{thm}\label{T3}
Let $I_{0}=\{(u^0_{3}(x),\,u^0_{4}(x))\mid 0\leq x\leq 1\}.$
\begin{enumerate}[\rm (i)]
\item Assume that there exists a point $(U_{0},\,V_{0})\in W^{s}$ such that $0<U_{0}<U_{s}$ and 
$I_{0}\subset (U_{0},\,\infty)\times (V_{0},\,\infty)$.
Then $(u(x,t),\,v(x,t))\to (u_{1},\,v_{1})$ uniformly on $0\leq x\leq 1$ as $t\to +\infty$.
\item
Assume that there exists a point $(U_{0},V_{0})\in W^{s}$ such that $0<U_{0}<U_{s}$ and $I_{0}\subset (0,U_{0})\times (0,V_{0})$.
Then $(u(x,t),\,v(x,t))\to (u_{0},v_{0})$ uniformly on $0\leq x\leq 1$ as $t\to +\infty$.
\end{enumerate}
\end{thm}

This theorem shows that  in order to obtain a nontrivial pattern, we must choose  the initial data  $I_{0}$ that are not contained in the rectangles stated in the theorem.


\subsubsection{Spatially nonhomogeneous stationary solutions}

Next, we construct a special stationary solution of \eqref{originalproblem} having a gradient-like form.

\begin{figure}[t]
 \begin{center}
   \includegraphics[width=120mm]{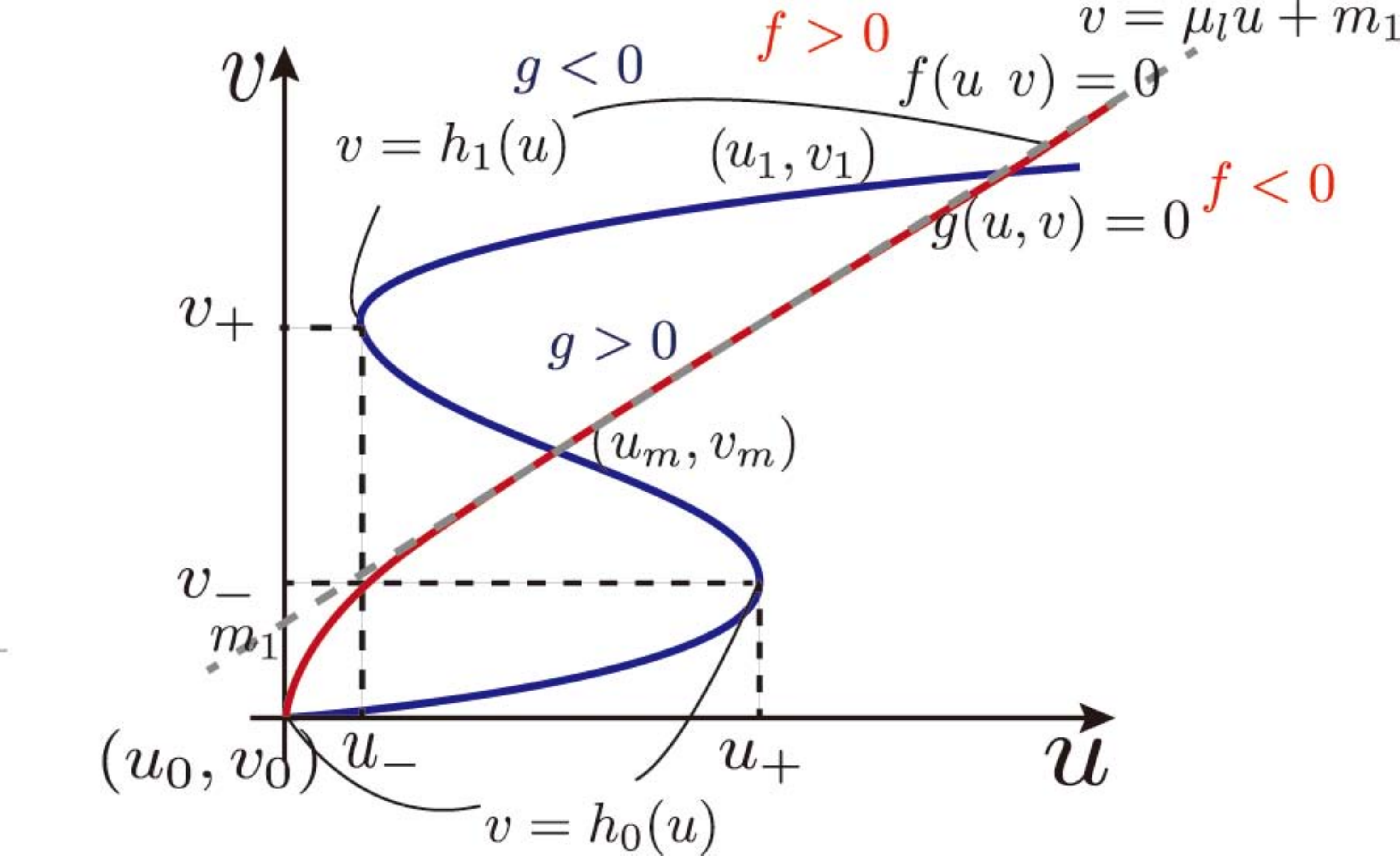}
\end{center}
\caption{Three constant solutions and two inflection points}
\label{zu2}
\end{figure}

Notice that $g(u,v)=0$ defines three smooth functions on the intervals $0 \leq u \leq u_{+}$,
$u_{-}\leq u\leq u_{+}$ and $u_{-}<u<+\infty$,
respectively (see Figure \ref{zu2}). Therefore, for each $u_{-}<u<u_{+}$,
the equation $g(u,v)=0$ has exactly three roots
$v=h_{0}(u),\,v=h_{m}(u)$ and $v=h_{1}(u)$ with
$h_{0}(u)<h_{m}(u)<h_{1}(u)$.
We can extend $h_{0}(u)$,\,$h_{m}(u)$ and $h_{1}(u)$ up to the end points of the intervals, so that $h_{0}(u_{+})=h_{m}(u_{+}),\,h_{m}(u_{-})=h_{1}(u_{-})$. It holds that
\begin{align}\label{h}
\left\{
\begin{aligned}
&v=h_{0}(u)&&\quad(u_{0}\leq u<u_{+},\,v_{0} \leq v<v_{-}),\\
&v=h_{m}(u)&&\quad(u_{-}< u < u_{+},\,  v_{-} < v< v_{+}),\\
&v=h_{1}(u)&&\quad(u_{-}<u<+\infty,\,v_{+}<v<+\infty).
\end{aligned}
\right.
\end{align}

\begin{definition}
A gradient-like pattern of system \eqref{originalproblem} is a stationary solution satisfying the following conditions:
 \begin{enumerate}[\rm (1)]
 \item $u$ is strictly monotone increasing;
 \item there exists an $l$ such that $0<l<1$ and 
 \begin{align*}
 v(x)=\left\{
 \begin{aligned}
 &h_{0}(u(x))&&\quad \text{\rm for}\,\, 0\leq x <l,\\
 &h_{1}(u(x))&&\quad \text{\rm for}\,\, l< x \leq1.\\
 \end{aligned}
 \right.
 \end{align*}
 \end{enumerate}
\end{definition}

Since $v$ has a jump discontinuity at $x=l$,
we require that $u(x)$ is of class $C^{1}([0,1])\cap C^{2}([0,l)\cup (l,1])$
and $u$ satisfies the second order differential equation in the sense of distribution.

Such a solution is important,
because we can construct other types of stationary solutions of \eqref{originalproblem}
by reflection,
periodic extension,
and rescaling 
(see Section \ref{SectionPattern}).

To construct a gradient-like solution we search for a continuously differentiable solution of the following boundary value problem with $F_{j}(u)=f(u,h_{j}(u))$
for $j=0,1$
\begin{align}\label{problem}
\left\{
\begin{aligned}
&\dfrac{d^{2}u}{dx^{2}}+\gamma F_{0}(u)=0 && \text{for}\,\, 0<x<l,\\
&\dfrac{d^{2}u}{dx^{2}}+\gamma F_{1}(u)=0 && \text{for}\,\, l<x<1,\\
&\dfrac{du}{dx}(0)=\dfrac{du}{dx}(1)=0,\,u(l)=\beta. &&
\end{aligned}
\right.
\end{align}

To state the result,
we define
\begin{align*}
\F_{0}(W) = \int_{u_{0}}^{W}F_{0}(z)dz & \qquad \text{for}\,\,u_{0}\leq W<u_{+}, \\
\F_{1}(W) =\int_{u_{1}}^{W}F_{1}(z)dz & \qquad \text{for}\,\, u_{-}< W \leq u_{1}.
\end{align*}
Note that $\F_{0}(W)<0$ for $ u_0<W<u_+$ and $\F_{1}(W)<0$ for $ u_- <W < u_1$.


\begin{thm}\label{MAIN}
Assume that $\mu_{3}$ satisfies \eqref{A3}.
Then,  for each  
$u_{-}<\beta<\min\{u_{+},u_{1}\}$ 
and
$0<m<\min\{\sqrt{2|\F_{0}(\beta)|},\sqrt{2|\F_{1}(\beta)|}\}$,
there exist a continuous function
$\gamma(\beta,m)>0,\,0<l(\beta,m)<1$,
and a continuously differentiable function
$u(x;\beta,m)$ on $0\leq x \leq 1$ such that $u=u(x;\beta,m)$  is a strictly  monotone increasing solution of \eqref{problem} for $\gamma=\gamma(\beta,m),\,l=l(\beta,m)$ and $m=\sqrt{\gamma(\beta,m)}u'(l(\beta,m);\beta,m)$.
Moreover,
\begin{enumerate}[\rm(i)]
\item\, $\gamma(\beta,m)\to 0$ and $u(x;\beta,m) \to\beta$ uniformly in $[0,1]$ as $m\downarrow 0$.
\item\,$\gamma(\beta,m)\to +\infty$ and $u'(l(\beta,m);\beta,m)\to+\infty$
as $m \uparrow \min\{\sqrt{2|\F_{0}(\beta)|},\sqrt{2|\F_{1}(\beta)|}\}$.
In addition,
\begin{enumerate}[\rm (a)]
\item if $\F_{0}(\beta)>\F_{1}(\beta)$,
then $l(\beta,m)\to\,1$ and $u(x;\beta,m)\to u_{0}$ locally uniformly in $[0,1)$,
\item if $\F_{0}(\beta)=\F_{1}(\beta)$,
then
$l(\beta,m)\to\,\dfrac{\sqrt{|F_{1}'(u_{1})|}}{\sqrt{|F_{1}'(u_{1})|}+\sqrt{|F_{0}'(u_{0})|}}=l^{*}$ and 
\begin{align*}
u(x;\beta,m)\to
\begin{cases}
u_{0}\quad\text{locally uniformly in}\ [0,l^{*}),\\
u_{1}\quad\text{locally uniformly in}\ (l^{*},1],
\end{cases}
\end{align*}
\item if $\F_{0}(\beta)<\F_{1}(\beta)$,
then $l(\beta,m)\to\,0$ and $u(x;\beta,m)\to u_{1}$ locally uniformly in $(0,1]$.
\end{enumerate}
\end{enumerate}
\end{thm}

Next,
we
study the uniqueness of gradient-like solution of \eqref{originalproblem} for a given $\gamma>0$.
Let
\begin{align*}
\triangle=\left|
\begin{array}{cc}
f_{u}(u,v)&f_{v}(u,v)\\
g_{u}(u,v)&g_{v}(u,v)
\end{array}
\right|.
\end{align*}
We can easily prove that there exist $\beta_{\triangle}$ and $B_{\triangle}$  satisfying 
$u_{0}<u_{-}<\beta_{\triangle}<B_{\triangle}<\min\{u_{+},u_{1}\}$  such that
if $\beta_{\triangle}<\beta<B_{\triangle}$  and $g(u,v)=0$,
then $\triangle>0$ .


\begin{thm}
{\rm (i)}
For each
$\gamma>0$,
if $\beta_{\triangle}<\beta<B_{\triangle}$, then boundary value problem \eqref{problem} has a unique monotone increasing solution.
{\rm (ii)}
On the other hand,
for each $\beta\in (u_{-},\,\min\{u_{+},\,u_{1}\})$,
there exists $\gamma^{*}>0$ such that \eqref{problem} has a unique monotone increasing solution whenever $\gamma>\gamma^{*}$. The same assertions hold true for monotone decreading solutions.
\end{thm}


\section{Initial-boundary value problem}\label{SectionExistence}

\noindent
In this section we prove Theorem \ref{ExistBounded}. First we apply the results on the existence of mild and classical solutions of the system of reaction-diffusion equations coupled with ordinary differential equations. 
They provide the local-in-time existence of solutions in the spaces of continuously differentiable and $\alpha$-H{\"o}lder functions. We refer to a generic form of the system of  $k$ reaction-diffusion equations coupled with $m$ ordinary differential equations given by \eqref{recmodel}, where $U=(u_{1}, u_{2} , \dots , u_{k})$, $V=(v_{1}, v_{2} , \dots , v_{m})$,  $F:\R^{N} \to \R^{k}$ and $G: \R^{N} \to \R^{m}$ with $N=k+m$ and  \begin{align*}
& F(U,V)=(f_{1}(U,V),\,\dots,\,f_{k}(U,V)),\\
& G(U,V)=(g_{1}(U,V),\,\dots,\,g_{m}(U,V)).
 \end{align*}
First, we recall that a {\it mild} solution  $(U,V)$  of problem \eqref{recmodel} on a time interval 
$[0,T)$   and with initial data $(U^0,V^0)\in \left( L^\infty (\Omega) \right)^N$
are  measurable functions $u_i,v_i :\Omega \times (0,T) \to \R$ satisfying the following system of
integral equations  
\begin{align}
& u_i(x,t)= S_i(t)u_{i}^0(x)+ \int_0^t S_i(t-s) (f_i(U(x,s),V(x,s))) \,ds  \quad \text{for}\ i=1,\dots, k, \\
&v_i(x,t)= v_{i}^0(x)+ \int_0^t g_i(U(x,s),V(x,s))\,ds \quad \text{for}\ i=1,\dots, m, 
\end{align}
where $S_i(t)$ is the semigroup of linear operators associated with the equation $z_t=d_i \Delta z -  \mu_i z$ %
in the domain  $\Omega$,  under the homogeneous Neumann boundary conditions.

\begin{prop} \label{proposition_existence}
Assume that $u_{i}^0, v_{i}^0 \in L^\infty (\Omega)$. %
Then, there exists $T > 0$ such
 that the initial-boundary value problem
 \eqref{recmodel} has a unique local-in-time mild
 solution $$u_i, v_i \in L^\infty \big([0, T],\, L^\infty (\Omega)
 \big).$$
If the initial data are more regular, i.e., $u_{i}^0 \in C^{2+\alpha}(\overline{\Omega})$,  $v_{i}^0 \in C^\alpha
 (\overline{\Omega})$ for some $\alpha \in(0,1)$
 and $\partial_\nu u_{i}^0 = 0$ on $\partial \Omega$, then the mild
 solution of problem is smooth and satisfies 
\begin{eqnarray*}
&&u_{i}\in C^{2+\alpha,1+\alpha/2} ( \overline{\Omega}\times [0,T]) \quad \text{for}\,\,\,\, i=1, \dots , k,\\
&&v_{i}\in C^{\alpha,1+\alpha/2} ( \overline{\Omega} \times [0,T]) \quad \text{for}\,\,\,\, i=1, \dots , m.
\end{eqnarray*}
\end{prop}
\noindent For the proof,  we refer to the lecture notes by Rothe \cite[Theorem 1, p.~111]{Rothe}, as well as to \cite[Theorem 2.1]{GSV} for studies of general
reaction-diffusion-ODE systems in the H\"older spaces.  
Our smoothness assumptions on the initial functions now guarantee the existence of a unique classical solution. (An elementary proof of Theorem \ref{ExistBounded} in the case of spatial dimension one is given in \cite[Appendix]{madokathesis}.)

Next,  we show the boundedness of solutions of the initial-boundary value problem \eqref{model} in a bounded domain $\Omega \subset \R^{n}$. We note that the nonnegativity of solutions for nonnegative initial conditions  is a consequence of the maximum principle. Adding the first two equations of  \eqref{model} and taking $\mu_{B}=\min\{\mu_{1},\mu_{2}\}$,
 we obtain
\begin{align}\label{315}
0< u_{1}(x,t)+u_{2}(x,t)\leq \dfrac{m_{1}}{\mu_{B}}+\max\{u^0_{1}(x),u^0_{2}(x)\}\quad \text{for}\,\,x\in\Omega\quad\text{and}\quad t>0.
\end{align}
Since the solutions are nonnegative, the above inequality implies that both $u_{1}$ and $u_{2}$ are uniformly bounded. 

Using this estimate, we show  the boundedness of $u_{3}(x,t)$ and $u_{4}(x,t)$. First,  let
\begin{align*}
&F(u_{3},u_{4};x,t)=-\mu_{3}u_{3}-b u_{3} u_{1}(x,t)+d
 u_{2}(x,t)+u_{4},\\
&G(u_{3},u_{4})=-\delta u_{4}+\dfrac{m_{4}u_{3}}{1+\sigma u_{4}^2-\beta_{l} u_{4}}.
\end{align*}
Note that $F(u_{3},u_{4})<0$ holds if and only if
 \begin{align*}
u_{4}&<\mu_{3}u_{3}+b u_{3} u_{1}(x,t)-du_{2}(x,t).
\end{align*}
Let $M= m_{1}/\mu_{B}+\max\{u^0_{1}(x),u^0_{2}(x)\}$.
Since $u_{2}\leq M$, the inequality
$$u_{4}<\mu_{3}u_{3}-dM$$
implies $F(u_{3},u_{4};t,x)<0$. 

\begin{figure}[t]
 \begin{center}
  \includegraphics[width=90mm]{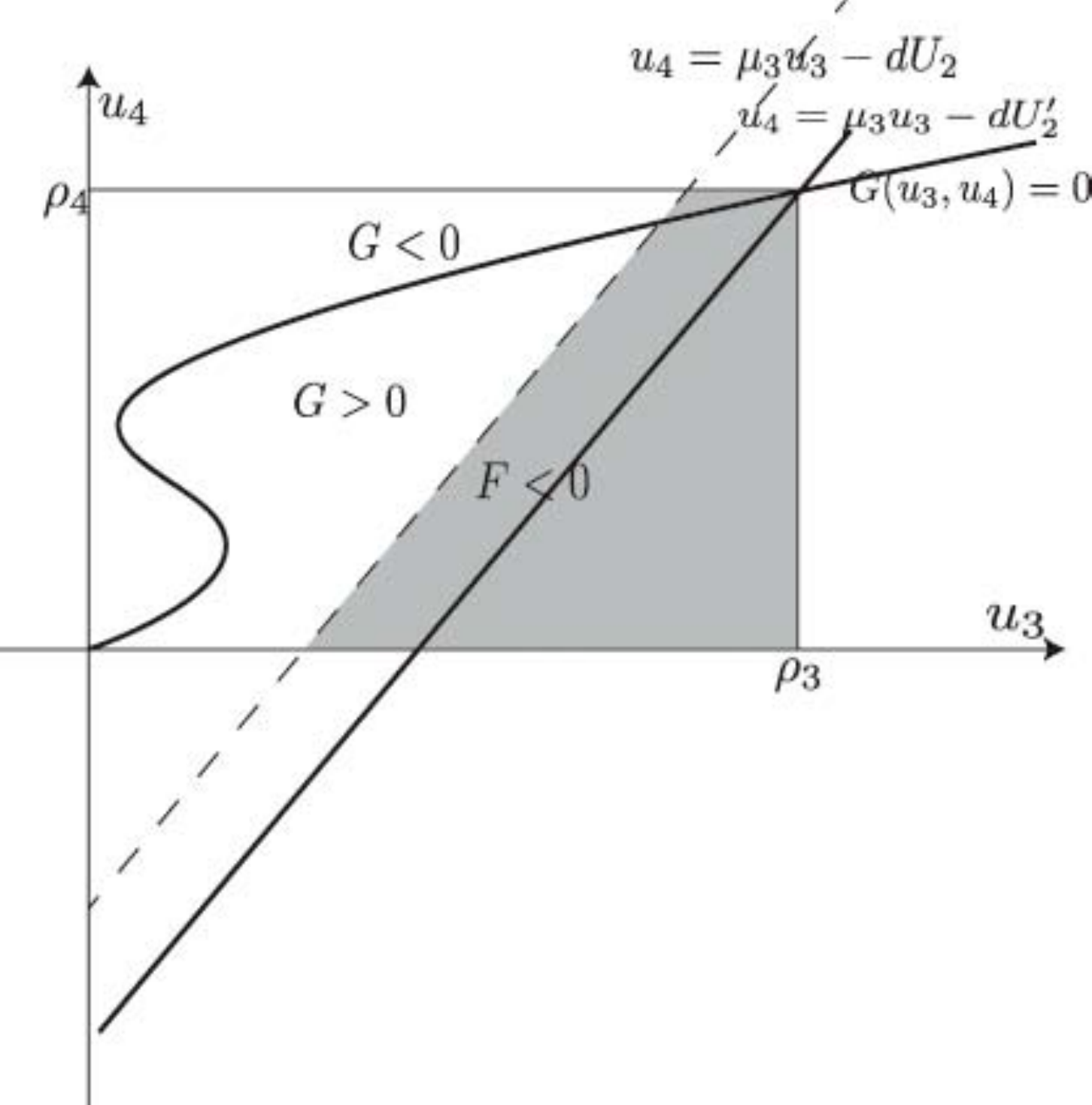}
 \end{center}
 \caption{Invariant rectangle for $(u_3,u_4)$}
 \label{rectangle}
\end{figure}

Clearly, we can find a positive constant $M'\geq M$ so large that the straight line $u_4=\mu_3 u_3 - dM^\prime$ and the curve $G(u_3,u_4)=0$ intersect at exactly one point $(\rho_*, \mu_3\rho_+-d M')$ in the first quadrant of the $u_3u_4$-plane  (see Figure \ref{rectangle}). Then, we obtain $(i)$ $F(\rho_3,u_4;x,t)<0$ for $0\leq u_4\leq \rho_4$, and $(ii)$ $G(u_3,\rho_4)>0$ for $0\leq u_3< \rho_4$. Note also that $\rho_4 \rightarrow +\infty$ as  $\rho_3 \rightarrow +\infty$, which provides arbitrarily large invariant sets.

Now, we take  $\rho_{3}>\rho_*$ so large that the rectangle ${\Sigma}=(0,\rho_{3})\times (0,\rho_{4})$ satisfies
$\{(u^0_{3}(x),u^0_{4}(x)) \mid x\in\overline{\Omega}\}\subset \Sigma$, where $\rho_{4}=\mu_{3}\rho_{3}-dM'$.
Since the vector field $\left( F(u_{3},u_{4};x,t), G(u_3,u_4) \right)$ points inside the set $\Sigma$, due to the maximum principle we obtain that
$0<u_{3}(x,t)<\rho_{3}$ and $0<u_{4}(x,t)<\rho_{4}$  for all $x\in \overline{\Omega}$. (See \cite{Smoller} for details of the framework of invariant rectangles for reaction-diffusion equations).

The proof of Theorem \ref{ExistBounded} is now complete.
\qed


\section{Initial values leading to uniform steady-states}\label{SectionKinetics}
\subsection{Comparison principle}
For the reduced two-component system \eqref{originalproblem} we have a comparison principle.
\begin{thm} \label{thm-comp}
If the two sets of initial data satisfy the inequalities $u_{1}^0(x) \ge u_{0}^0(x)$ and 
 $v_{1}^0(x) \ge v_{0}^0(x)$ for all $x\in \overline \Omega$, then 
$$u_{1}(x,t) \ge u_{0}(x,t) \quad \text{and} \quad  
 v_{1}(x,t) \ge v_{0}(x,t) \qquad 
 \text{ for all x } \in \overline\Omega, \ t\geq0.
 $$ 
\end{thm}
\begin{prf}
Let $\phi$ and $\psi$ be defined by
$$\phi(x,t)=u_1(x,t)-u_0(x,t), \quad \psi(x,t)=v_1(x,t)-v_0(x,t).$$
Then, by the mean value theorem we see that $\phi$ and $\psi$ satisfy
\begin{align*}
\left\{
\begin{aligned}
&\dfrac{\partial \phi}{\partial t} =\dfrac{1}{\gamma}\Delta \phi - a(x,t)\phi + \psi && \text{for}\ x\in \Omega,\ t >0,\\
&\dfrac{\partial \psi}{\partial t}=b(x,t)\phi + c(x,t)\psi &&\text{for}\ x\in \Omega,\ t >0,\\
&\dfrac{\partial \phi}{\partial \nu}=0 &&\text{for}\ x \in \partial \Omega,\ t >0,
\end{aligned}
\right.
\end{align*}
together with initial conditions
\begin{align*}
&\phi(x,0)=u_1^0(x)-u_0^0(x),\  \psi(x,0)=v^0_{1}(x)-v_0^0(x) \quad \text{for}\ x \in \Omega,
\end{align*}
where
\begin{align*}
a(x,t)&= \int_0^1 f_u\left(u_0(x,t)+s\phi(x,t),v_0(x,t)+s\psi(x,t)\right)ds,\\
b(x,t)&= \int_0^1 g_u\left(u_0(x,t)+s\phi(x,t),v_0(x,t)+s\psi(x,t)\right)ds,\\
c(x,t)&= \int_0^1 g_v\left(u_0(x,t)+s\phi(x,t),v_0(x,t)+s\psi(x,t)\right)ds.
\end{align*}
Moreover, it follows that 
$$a(x,t)>0, \quad b(x,t)>0.$$
(The sign of $c(x,t)$ is not definite.) 
The latter conditions together with the maximum principle yield the nonnegativity of $\phi(x,t)$ and $\psi(x,t)$ for nonnegative initial conditions $\phi(x,0)$ and $\psi(x,0)$. This completes the proof of Theorem \ref{thm-comp}.
\end{prf}


\subsection{Stable manifold for the kinetic system}

\begin{figure}[t]
 \begin{center}
   \includegraphics[width=80mm]{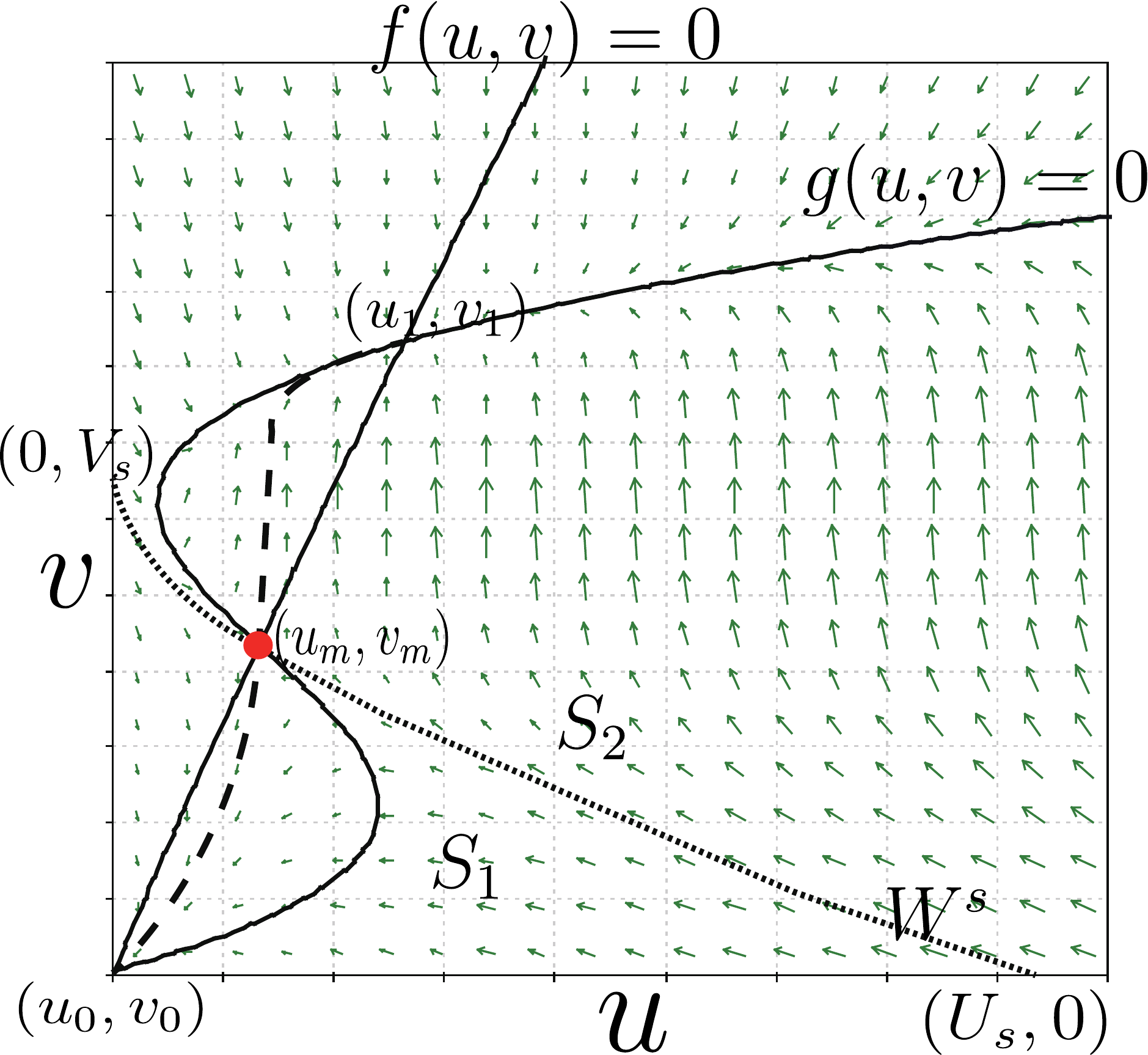}
 \end{center}
\caption{Vector field and the stable manifold $W^{s}$}
 \label{pplane}
\end{figure}

We consider kinetic system \eqref{kineticsystem} and recall that under assumptions \eqref{A2},\,\eqref{A3},
the system has three equilibria $(u_{0},v_{0}),\,(u_{m},v_{m})$ and $(u_{1},v_{1})$, such that  $u_{0}=0<u_{m}<u_{1}$ and $v_{0}=0<v_{m}<v_{1}$. The Jacobi matrix $A$ evaluated at $(u_m,v_m)$
\begin{align*}
A=\left(
\begin{array}{cc}
f_{u}(u_{m},v_{m}) &f_{v}(u_{m},v_{m})\\
g_{u}(u_{m},v_{m}) &g_{v}(u_{m},v_{m})
\end{array}
\right)
\end{align*}
has two real eigenvalues 
$\mu < 0 < \lambda$ 
and the equilibrium 
$(u_{m},
v_{m})$ 
is a saddle point. 
It is verified by straighforward calculations that
\begin{align*}
&f_{u}(u_{m},v_{m})=-\mu_{3}-\dfrac{\mu_{1}\mu_{2}bm_{1}(\mu_{2}+d)}{(\mu_{1}(\mu_{2}+d)+\mu_{2}bu_{m})^{2}}<0,\\
&f_{v}(u_{m},v_{m})=1>0,\\
&g_{u}(u_{m},v_{m})=\dfrac{m_{2}}{a+\sigma v_{m}^{2}-\beta_{l} v_{m}}>0,\\
&g_{v}(u_{m},v_{m})=-\dfrac{\delta(1+3\sigma v_{m}^{2}-2\beta_{l} v_{m})}{1+\sigma v_{m}-\beta_{l} v_{m}}>0.
\end{align*}
Hence, 
the equilibrium 
$(u_{m},
v_{m})$ of system  \eqref{kineticsystem} 
has a stable manifold $W^{s}$ and an
unstable manifold $W^{u}$.
(See Figure \ref{pplane}.)

Now, we are ready to prove Lemma \ref{saddle} which is a premise of Theorem \ref{T3}.

\noindent
{\bf Proof of Lemma \ref{saddle}.}
Let $(u_{\mu}, v_{\mu})$ 
be an eigenvector of $A$ belonging to $\mu$. 
Then,
\begin{align} \label{5}
\left\{
\begin{aligned}
& f_{u}(u_{m},
v_{m}) u_{\mu}
+f_{v}(u_{m},
 v_{m}) v_{\mu} = \mu u_{\mu}, \\
& g_{u}(u_{m},
v_{m}) u_{\mu}
+g_{v}(u_{m},
 v_{m}) v_{\mu} = \mu v_{\mu}.
\end{aligned}
\right.
\end{align}
Note that 
$u_{\mu} v_{\mu} < 0$
 since it holds that (i) 
$g_{u}(u_{m}, v_{m}) u_{\mu} v_{\mu}
 = (\mu - g_{v}(u_{m},v_{m}))(v_{\mu})^{2}$ 
and
(ii) 
$g_{u}(u_{m}, v_{m}) > 0$, 
$\mu - g_{v}(u_{m}, v_{m}) < 0$. 
Hence, 
choosing
$u_{\mu} = 1$ 
yields $v_{\mu} < 0$.
From \eqref{5} 
we observe that the angle between this eigenvector 
$(u_{\mu},v_{\mu})$ 
and the normal vector 
$( f_{u}(u_{m}, v_{m}), f_{v}(u_{m}, v_{m}))$ 
to 
$f (u,v) = 0$ is greater than $\pi / 2$, 
while the angle between $(u_{\mu},
v_{\mu})$
and the normal vector 
$(g_{u}(u_{m}, v_{m}), g_{v}(u_{m}, v_{m}))$ 
to the curve $g(u, v) = 0$ 
is smaller than $\pi / 2$.
For a positive constant $\rho$, 
define a disk $D_{\rho} = \{ (u, v) \in \R^{2} \mid j
 (u-u_{m})^{2} +(v-v_{m})^{2} < \rho^{2} \}$.
If $\rho$ is sufficiently small then 
$D_{\rho}$ is divided into four subdomains $R_{1}$, 
$R_{2}$, 
$R_{3}$, 
$R_{4}$ and the portions of $f (u, v) = 0$ 
and $g(u, v) = 0$, 
where
\begin{align*}
& R_{1} = \{(u, v) \mid f (u, v) < 0, g(u, v) > 0 \},
\quad 
R_{2} =  \{(u, v) \mid f (u, v) > 0, g(u, v) > 0 \},\\
& R_{3} = \{(u, v) \mid f (u, v) > 0, g(u, v) < 0 \},
\quad
 R_{4} =  \{(u, v) \mid f (u, v) < 0, g(u, v)< 0 \}.
\end{align*}
The consideration above implies that the tangent 
to the stable manifold $W^{s}$ for the equilibrium
$(u_{m}, v_{m})$ 
is included in 
$R_{1} \cup \{ (u_{m}, v_{m}) \} \cup R_{3}$, 
that is, 
$W^{s} \cap\ (D_{\rho} \setminus \{(u_{m}, v_{m}) \} )$ 
$\subset (R_{1} \cup R_{3})$.
First we consider the part $W^{s} \cap R_{1} \cap D_{\rho}$.
It is convenient to reverse the orientation of
time variable $t \mapsto -t$ and put 
$\tilde{u}(t) = u(-t)$, 
$\tilde{v}(t) = v(-t)$. 
Then,
\begin{align} \label{6}
\left\{
\begin{aligned}
& \frac{d \tilde{u}}{dt}
 =
-f ( \tilde{u},
\tilde{v}), \\
& \frac{d \tilde{v}}{dt}
 =
-g ( \tilde{u},
\tilde{v}).
\end{aligned}
\right.
\end{align}
We choose an initial data $(\tilde u(0), \tilde v(0) )= (u_{\ast}, v_{\ast}) \in W^{s} \cap R_{1} \cap
D_{\rho}$.
Since $\tilde{u}'(t) > 0$ 
and $\tilde{v}'(t) < 0$ in $R_{1}$, 
there are three possibilities:
\begin{enumerate}
\item[(A)] \,There exists a positive constant $T > 0$ 
such that 
$( \tilde{u}(t), \tilde{v}(t)) \in R_{1}$ for $0 \le t < T$ and
$( \tilde{u}(T), \tilde{v}(T))$ 
is on the curve $g(u,v) = 0$, 
so that $\tilde{v}'(T) = 0$. 
\item[(B)]\,There exists a positive constant $T > 0$ 
such that 
$( \tilde{u}(t), \tilde{v}(t)) \in R_{1}$ for $0 \le t < T$ and
$\tilde{u}(T) > 0$, 
$\tilde{v}(T) = 0$.
\item[(C)]\, $( \tilde{u}(t), \tilde{v}(t)) \in R_{1} \cap Q_{1}$
 for all $t$ in the maximum existence interval $0 \le t <T_{M}$ and
$\tilde{u}(t) \to +\infty$ as $t \to T_{M}$.
\end{enumerate}
On the boundary between $R_{1}$ 
and $R_{4}$, the vector field $(-f ( \tilde{u}, \tilde{v}),
-g( \tilde{u}, \tilde{v}))$ points to the interior of $R_{1}$.
Therefore, 
the orbit in $R_{1}$ cannot enter $R_{4}$ through the boundary. 
Hence,
even if (A) occurs, the solution $( \tilde{u}(t),
\tilde{v}(t))$ 
returns to $R_{1}$. 
Therefore, 
we assume case (C) and derive a contradiction. 
From
\begin{align*}
0 <
\frac{d \tilde{u}}{
dt} 
= -
f ( \tilde{u}, \tilde{v}) 
= -\tilde{v} + \mu_{3} \tilde{u} 
+ m_{1} 
-\frac{ \mu_{1}( \mu_{2} + d)m_{1}}
{\mu_{1}(\mu_{2} + d) + \mu_{2} b \tilde{u} }
\le \mu_{3} \tilde{u} + m_{1}
\end{align*}
we obtain
\begin{align*}
u_{\ast} \le \tilde{u}(t) \le e^{\mu_{3} t} y_{\ast} +
\frac{m_{1}}{\mu_{3}}
\left(
e^{\mu_{3} t}
-1
\right).
\end{align*}
This leads to 
$T_{M} = \infty$ in case (C). 
Since $0 < \tilde{v}(t) \le v_{\ast}$ 
for all $0 \le t < 1$, 
it holds that for $t > 0$ 
\begin{align*}
\frac{d \tilde{v} }{dt}
(t) = 
-g( \tilde{u}(t), \tilde{v}(t)) 
= \delta \tilde{v}(t) 
-\frac{ m_{2} \tilde{u}(t)}
{
1 + \sigma \tilde{v}(t)^{2} - \beta_{l} \tilde{v}(t) 
}
\le
\delta v_{\ast} - \frac{m_{2}}
{1 + \sigma v_{\ast}^{2} }
\tilde{u}(t).
\end{align*}
Then, 
choosing $t_{0}$ so large that 
$\delta v_{\ast}- m_{2} \tilde{u}(t_{0})/(1 + \sigma v_{\ast}^ 2) 
< -\delta v_{\ast}$,
by $ u^\prime(t)>0$ we obtain
\begin{align*}
0 < \tilde{v}(t) \le \tilde{v}(t_{0}) - \delta v_{\ast} (t- t_{0})
\qquad \text{for} \quad t > t_{0},
\end{align*}
which yields a contradiction if 
$t > t_{0} + \tilde{v}(t_{0})/(\delta v_{\ast})$. 
We conclude that case (B) holds.

Furthermore, 
we can show that solution 
$( \tilde{u}(t), \tilde{v}(t))$ of \eqref{6} with initial data
$(u_{\ast},v_{\ast})$ $\in R_{3} \cap D_{\rho}$ 
reaches the positive $v$-axis in finite time by using the estimates
\begin{align}
 0 >
\frac{d \tilde{u}}{
dt} 
& = -f (\tilde{u}, \tilde{v}) \ge -\tilde{v}- m_{1}, \label{7} \\ 
\max\{ 0, \delta v_{\ast} \} -\frac{m_{2}}{c_{0}} \tilde{u} 
& < 
\frac{d \tilde{v}}{dt}
= -g(\tilde{u}, \tilde{v}) \le \delta \tilde{v}, \label{8}
\end{align}
where $c_{0} = \min_{v} (1 + \sigma v^{2} - \beta_{l} v)> 0$. 
Indeed, 
\eqref{8} implies that the maximum existence interval
for solutions which stay in 
$R_{3}$ is equal to $0 \le t < \infty$. 
If $\tilde{v}(t) \to +\infty$ as $t \to \infty$, 
then \eqref{7} yields that 
$\tilde{u}(t)$ cannot remain positive for all $t > 0$. 
Thus in this case 
$\tilde{u}(T_{0}) = 0$ for some $0 < T_{0} < \infty$, 
which in turn implies $v_{\ast} < \tilde{v}(T_{0}) < \infty$. 
If $\tilde{v}(t) \le K$ for all $t > 0$, we obtain  by \eqref{7}
$$
\frac{d \tilde{u}}{dt} \ge 
 -K - m_{2} \frac{\tilde{u}(t)}
{c_{0} } > -\frac{K}{2} 
$$
provided that $t$ is sufficiently large. 
Hence,
$u(T_{1}) = 0$ for some $0 < T_{1} < \infty$. 
Consequently we see that there exists a finite $T$ such that
$\tilde{u}(T) = 0$ and $0 < \tilde{v}(T) < \infty$. 
\qed


\begin{lem}\label{kinetic}
Let ${\cal{Q}}=S_{1}\cup W^{s}\cup S_{2}$,
where $S_{1}$ and $S_{2}$ are disjoint open sets and $(0,0)\in\bar{S_{1}}$.
Let $(u(t),v(t))$ be a solution of the kinetic system  \eqref{kineticsystem}.
Then
\begin{enumerate}[\rm (a)]
\item if $(u(0),v(0))\in S_{1}$,\,then $(u(t),v(t))\to (0,0)$ as $t\to +\infty$;
\item if $(u(0),v(0))\in S_{2}$,\,then $(u(t),v(t))\to (u_{1},v_{1})$ as $t\to +\infty$.
\end{enumerate}
\end{lem}

%
\begin{prf}
Solutions of the kinetic system  \eqref{kineticsystem} remain bounded as $t\to \infty$ and they cannot be periodic due to the structure of the vector field $(f(u,v),g(u,v))$. Indeed,
the nullclines  divide the first quadrant $\cQ$ into four open sets
\begin{align*}
{\cR}_{1}=\{(u,v) \mid f(u,v)<0,\,g(u,v)>0\},\quad {\cR}_{2}=\{(u,v) \mid f(u,v)<0,\,g(u,v)<0\},\\
{\cR}_{3}=\{(u,v) \mid f(u,v)>0,\,g(u,v)<0\},\quad {\cR}_{4}=\{(u,v) \mid f(u,v)>0,\,g(u,v)>0\},
\end{align*}
where ${\cR}_{1},\,{\cR}_{3}$ and ${\cR}_{4}$ are connected,
but ${\cR}_{2}$ has two connected components ${\cR}_{2,0}$ and ${\cR}_{2,1}$
such that $(0,0)\in\overline{{\cR}_{2,0}}$, whereas $(u_{1},v_{1})\in\overline{{\cR}_{2,1}}$.
We observe that ${\cR}_{2,0}$,\,${\cR}_{2,1}$ and ${\cR}_{4}$ are invariant sets,
moreover if $(u(0),v(0))\in {\cR}_{2,0}$ then $(u(t),v(t))\to (0,0)$ as $t\to +\infty$;
if $(u(0),v(0))\in {\cR}_{2,1}\cup {\cR}_{4}$ then $(u(t),v(t))\to (u_{1},v_{1})$ as $t\to \infty$.
If $(u(0),v(0))\in {\cR}_{1}\cup {\cR}_{3}$,
then $(u(t),v(t))$ either leaves ${\cR}_{1}\cup {\cR}_{3}$ in finite time or stays there for all $t>0$
and converges to $(u_{m},v_{m})$ if $(u(0),v(0))\in{\cR}_{1}$ or to $(0,0)$ if $(u(0),v(0))\in {\cR}_{3}$ as $t\to +\infty$.
Now, we consider Case (a).
Note that
$$S_{1}\backslash(\{f=0\}\cup\{g=0\})=({\cR}_{1}\cap S_{1})\cup {\cR}_{2,0}\cup ({\cR}_{3}\cap S_{1}).$$
Since only solutions with $(u(0),v(0))\in W^{s}$ converge to $(u_{m},v_{m})$,
we conclude that for $(u(0),v(0))\in ({\cR}_{1}\cap S_{1})\cup({\cR}_{3}\cap S_{1})\cup {\cR}_{2,0}$,
$(u(t),v(t))\to (0,0)$ as $t\to+\infty$.

Case (b) is treated in the same way, and hence we omit the details.
\end{prf}


\subsection{Convergence to constant stationary solutions}
In this subsection, we prove Theorem \ref{T3}.
Case (i):
There exists a point $(U_{0},V_{0})\in W^{s}$ such that $I_{0}\subset (U_{0},\infty)\times(V_{0},\infty)$.
Then there are positive constants $\delta$ and $M$ with $\delta<M$ such that $I_{0}\subset (U_{0}+\delta,U_{0}+M)\times(V_{0}+\delta,V_{0}+M)$.
Let $(u^{1}(t),v^{1}(t))$ and  $(u^{2}(t),v^{2}(t))$ be solutions of kinetic system   
 \eqref{kineticsystem} for initial values 
$$(u^{1}(0),v^{1}(0))=(U_{0}+\delta,\,V_{0}+\delta)\quad \text{and}\quad (u^{2}(0),\,v^{2}(0))=(U_{0}+M,\,V_{0}+M),$$
respectively.
Then, $(u^{1}(t),\,v^{1}(t))$ is a solution of the initial-boundary value problem \eqref{originalproblem} with initial value $u^{1,0}(x)=U_{0}+\delta$,
$v^{1,0}(x)=V_{0}+\delta$, while $(u^{2}(t),v^{2}(t))$ is a solution of
\eqref{originalproblem} with initial value $u^{2,0}(x)=U_{0}+M,\,v^{2,0}(x)=V_{0}+M$.
Applying Theorem \ref{thm-comp}, we obtain that for all  $x\in[0,1]$ and $t\geq 0$,
\begin{align*}
u^{1}(t)\leq u(x,t)\leq u^{2}(t)\quad \text{and}\quad
v^{1}(t)\leq v(x,t)\leq v^{2}(t).
\end{align*}
Now we know that  $(u^{1}(0),v^{1}(0))$ and $(u^{2}(0),v^{2}(0))\in S_{2}$. Then Lemma \ref{kinetic} yields that
$(u^{1}(t),v^{1}(t))\to (u_1,v_1)$ and $(u^{2}(t),v^{2}(t))\to (u_1,v_1)$
as $t\to +\infty$.
Therefore, 
$(u(x,t),v(x,t)) \to (u_1,v_1)$ as $t\to +\infty$ uniformly on $\overline{\Omega}$.

Case (ii) is treated in the same way.
\qed


\section{Nonconstant stationary solutions}\label{SectionPattern}
\subsection{Preliminaries}
In this section
we construct stationary solutions of the two-equation problem \eqref{originalproblem} such that $u$ is monotone (increasing or decreasing). If $(u(x),v(x)) $ is a solution of \eqref{originalproblem} with $u(x)$ being monotone increasing, then $(u(1-x), v(1-x)) $ gives rise to a monotone decreasing solution. Therefore, we concentrate on monotone increasing solutions. Let us consider the boundary value problem
\begin{align} \label{beforechange}
\left\{
\begin{aligned}
\dfrac{1}{\gamma}\dfrac{d^{2}u}{dx^{2}}+f(u,v)&=0 &&\quad \text{for}\quad  0<x<1,\\
g(u,v)&=0  &&\quad \text{for}\quad  0<x<1,\\
\dfrac{du}{dx}&=0 &&\quad \text{at}\quad\,\,  x=0,1.
\end{aligned}
\right.
\end{align}
where 
$f(u,v)$ and $g(u,v)$ are defined by \eqref{2eqnonlinearities}.

Fix a $\beta$ in the interval $u_- < \beta < \min \{ u_+, u_1 \} $. We solve the second equation of \eqref{beforechange} for $v$ by switching between two branches $ v=h_0(u) $ and $ v=h_1(u) $ at $u=\beta$, i.e., $ v=h_0(u) $ if $ u < \beta $ and $ v=h_1(u) $ if $ u>\beta $. Hence, we define a function $ F_\beta (u) $ by
$$
F_\beta(u) = 
\begin{cases}
 F_0(u) & \text{if} \ u<\beta ,\\ F_1(u) & \text{if} \ u>\beta. \\ \end{cases}
$$
Now \eqref{beforechange} is reduced to the following boundary value problem for $u$ alone:
\begin{align} \label{5-2}
\left\{
\begin{aligned}
\dfrac{1}{\gamma}\dfrac{d^{2}u}{dx^{2}}+F_\beta(u)&=0 &&\quad \text{for}\quad  0<x<1,\\
\dfrac{du}{dx}&=0 &&\quad \text{at}\quad  x=0,1.
\end{aligned}
\right.
\end{align}
Since the nonlinear term $ F_\beta(u) $ has a jump discontinuity at $u=\beta$, we cannot expect a classical solution of \eqref{5-2}. In what follows, we always require a solution of \eqref{5-2} to be continuously differentiable on the interval $ 0 \leq x \leq 1$ and satisfies the first equation in the sense of distribution. The goal of this section is to prove that for any $\beta $ in the interval $ u_- < \beta < \min \{ u_+, u_1 \} $, there is such a solution. As a matter of fact, similar problems have been considered in order to construct a first approximate solution (or, outer solution) of reaction-diffusion systems where both species diffuse. For instance, see a classical paper by Mimura, Tabata and Hosono \cite{MTH}. Our problem does not satisfy all of the assumptions made in \cite{MTH} and hence we take a somewhat different approach applicable under weaker assumptions than theirs. We regard the coefficient $\gamma $ also as an unknown, and find a one-parameter family of solutions $(u(x;m), \gamma(m))$ for some parameter $m \in I $ where $I$ is an interval. It turns out that $\{ \gamma(m) \mid m \in I \} = (0,+\infty)$; and therefore, given $\gamma>0$, one can always find at least one monotone increasing solution of \eqref{5-2}.

\noindent

Assuming that $u(x)$ is strictly increasing,
we have an $l$ in the interval $(0,1)$ such that
\begin{align}\label{M}
\left\{
\begin{aligned}
&\dfrac{d^{2}u}{dx^{2}}+\gamma F_{0}(u)=0 && \text{for}\,\, 0<x<l,\\
&\dfrac{d^{2}u}{dx^{2}}+\gamma F_{1}(u)=0 && \text{for}\,\, l<x<1,\\
&\dfrac{du}{dx}(0)=\dfrac{du}{dx}(1)=0, &&\\
&u(l)=\beta. &&
\end{aligned}
\right.
\end{align}
Let us explain the procedure to obtain a one-parameter family of solutions of \eqref{5-2}.
\begin{enumerate}
\item[(i)]\,First, we scale the spatial variable $x$ by $x=l+y/\sqrt{\gamma}$ and 
define
$u(x)=U(y)$. 
By this change of variable, 
the interval $0 \leq x \leq 1$ is transformed onto the interval $-\sqrt{\gamma}l \leq y \leq \sqrt{\gamma}(1-l)$.
\item[(ii)]\,Problem \eqref{M} is now converted to the equivalent problem
\begin{align}\label{main}
\left\{
\begin{aligned}
&\dfrac{d^{2}U}{dy^{2}}+F_{0}(U)=0 \qquad\text{for}\,\, -\sqrt{\gamma}l < y < 0, \\ 
&\dfrac{d^{2}U}{dy^{2}}+F_{1}(U)=0 \qquad\text{for}\,\,  0 < y < \sqrt{\gamma}(1-l),\\ 
&\dfrac{dU}{dy}(-\sqrt{\gamma}l)=\dfrac{dU}{dy}(\sqrt{\gamma}(1-l))=0, \\ 
&U(0)=\beta \quad \text{and} \quad U'(y)>0 \quad \text{for}\,\,-\sqrt{\gamma}l < y < \sqrt{\gamma}(1-l).  
\end{aligned}
\right.
\end{align}
\item[(iii)]\,Instead of solving the boundary value problem \eqref{main},
we solve the following two initial value problems
\begin{align}\label{440}
\left\{
\begin{aligned}
&U''=-F_{0}(U)\quad\text{for}\,\, y<0,\\
&U(0)=\beta,\quad U'(0)=m
\end{aligned}
\right.
\end{align}
and
\begin{align}\label{441}
\left\{
\begin{aligned}
&U''=-F_{1}(U)\quad\text{for}\,\, y>0,\\
&U(0)=\beta,\quad U'(0)=m,
\end{aligned}
\right.
\end{align}
where $m$ is a given positive number.
\item[(iv)]\,Let $U_{0}$ be the solution of \eqref{440} and $U_{1}$ the solution of \eqref{441}.
If there exist  two positive constants $M$ and $N$ such that
\begin{enumerate}
\item[1)]\, $U_{0}'(y)>0$ for $-M<y<0$ and $U_{0}'(-M)=0$,
\item[2)]\, $U_{1}'(y)>0$ for $0\leq y<N$ and $U_{1}'(N)=0$,
\end{enumerate} 
then we define
\begin{align}\label{442}
l=\dfrac{M}{M+N}\quad \text{and} \quad \gamma=(M+N)^{2}.
\end{align}
Defining $U(y)$ by
\begin{align}
U(y)=\begin{cases}
U_{0}(y)\quad \text{for}\,\, -M\leq y\leq 0,\\
U_{1}(y)\quad \text{for}\,\, 0<y \leq N,
\end{cases}
\end{align}
we obtain a solution of \eqref{main} with $l$ and $\gamma$ determined by \eqref{442}.
Note that $ U \in C^1([-\sqrt{\gamma}\,l, \sqrt{\gamma}(1-l)] ) \cap C^2([-\sqrt{\gamma} l, 0]) \cap C^2([0, \sqrt{\gamma} (1-l) ])$, since $ F_0(U) $ and $ F_1(U) $ are both smooth up to $ u=\beta $.
\end{enumerate}  

Therefore, for the proof of Theorem \ref{MAIN}, we have to show that there do exist $M$ and $N$ satisfying 1) and 2) above, which is done in the following

\begin{lem} Given 
$\beta \in (u_{-}, \min\{u_{+},u_{1}\})$, let $m $ be any number in the interval $0 < m < \min \left\{ \sqrt{-2 \F_0(\beta) }, \sqrt{ -2 \F_1(\beta) } \right\}$. Then there exists a unique $k \in (0, \beta) $ and $ p \in (\beta, u_1) $ such that
\begin{align}\label{houteishiki}
\F_{0}(k)=
\dfrac{1}{2}m^{2}+\F_{0}(\beta),
\quad \F_{1}(p)=\dfrac{1}{2}m^{2}+\F_{1}(\beta).
\end{align}
Let
\begin{equation*}
M=\dfrac{1}{\sqrt{2}}\int_{k}^{\beta}\dfrac{dw}{\sqrt{\F_{0}(k)-\F_{0}(w)}},\quad
N=\dfrac{1}{\sqrt{2}}\int_{\beta}^{p}\dfrac{dw}{\sqrt{\F_{1}(p)-\F_{1}(w)}}.
\end{equation*}
Then the solution $ U_0(y) $ of \eqref{440} satisfies {\rm 1)} of\/ {\rm Step (iv)} with $k=U_0(-M)$, and the solution $ U_1(y) $ of \eqref{441} satisfies {\rm 2)} with $p=U_1(N)$.
\end{lem}


\begin{prf} 
We consider $ U_0 $ only, since $ U_1$ is treated in exactly the same way.
We multiply both sides of  \eqref{440} by $ U^\prime $ to obtain
\begin{align*}
\dfrac{dU}{dy} 
\dfrac{d^{2}U}{dy^{2}} 
+F_{0}(U)\dfrac{dU}{dy}=0.
\end{align*}
Recalling the definition of $ \F_0(U)$, we write this in the following form:
\begin{align}\label{dainyu}
 \dfrac{1}{2}\dfrac{d}{dy}\left(\dfrac{dU}{dy}\right)^{2}+\dfrac{d}{dy}\F_{0}(U)=0.
\end{align}
Now
we integrate both sides of \eqref{dainyu} over the interval $[y, 0]$ and obtain
\begin{align} \label{5-101}
\dfrac{1}{2}m^{2}+\F_{0}(\beta)
-\left\{\dfrac{1}{2}\left(\dfrac{dU}{dy}(y)\right)^{2}+\F_{0}(U(y))\right\}&=0,
\end{align} 
where we used the initial conditions $ U(0)=\beta$ and $ U^\prime(0)=m $.
Hence, a monotone increasing solution satisfies
$$
\dfrac{dU}{dy}(y)=\sqrt{m^{2}+2\{\F_{0}(\beta)-\F_{0}(U(y))\}}.
$$
This solution is well-defined as long as
$ E(U(y)) = m^{2}+2\{\F_{0}(\beta)-\F_{0}(U(y))$ is nonnegative.
Since $ E(U(0)) = m^2 > 0 $ and $E^\prime(U)=-\F'_{0}(U)>0$, as $y$ decreases from $0$, $E(U(y)) $ remains positive for a while, but decreases with the decrease of $y$ (note that $U(y)$ also decreases). Clearly, $U(y)$ is decreasing as $y$ decreases until it reaches the value $k$ for which
\begin{align}\label{beta}
m^{2}+2\{\F_{0}(\beta)-\F_{0}(k)\}=0
\end{align}
holds.
Since $ \F_0^\prime(U) < 0$, this algebraic relation determines $k$ uniquely in the interval $0<k<\beta$, and $k$ is a continuously differentiable function of $(\beta,m)$.
From \eqref{beta}, we notice that \eqref{5-101} can be written in the following form:
\begin{align*}
\dfrac{dU}{dy}(y)=\sqrt{2}\sqrt{\F_{0}(k)-\F_{0}(U(y))}.
\end{align*}
We integrate this equation to obtain $U(y)$ as the inverse function of
\begin{align}
\frac{1}{ \sqrt{2} } \int_u^\beta \frac{ dw } { \sqrt{ \F_0(k) - \F_0(w) } } = -y.
\end{align}
Observe that, as $ U \downarrow k$, the integral on the left-hand side of the identity above is convergent, since $ \F_0(w)=\F_0(k)+(F_0(k)+o(1))(w-k) $ and $ F_0(k) < 0$.
Therefore, we can define $M$ by
\begin{equation*}
M=M(\beta,m)=\dfrac{1}{\sqrt{2}}\int_{k}^{\beta}\dfrac{dw}{\sqrt{\F_{0}(k)-\F_{0}(w)}}.
\end{equation*}
(The dependence of $M$ on $m$ is by way of $k=k(\beta,m)$.)
This shows that $ U(-M)=k $ and $U^\prime(-M)=0$, which completes the proof of the lemma.
\end{prf}

\medskip

To summarize, for each $\beta \in (u_-, \min \{u_+, u_1 \} )$, we have constructed a one-parameter family of solutions $\left \{ u(x;\beta,m) \bigm | 0< m < \min \{ \sqrt{2} \F_0(\beta), \sqrt{ 2| \F_1(\beta) |} \} \right\} $ of \eqref{M} with
$$
l=l(\beta,m) = M(\beta,m)/(M(\beta,m)+N(\beta,m)), \ 
\gamma = \gamma(\beta,m) =(M(\beta,m)+N(\beta,m))^2.
$$
It is to be emphasized that $ M(\beta,m) $ and $ N(\beta,m) $ are continuous in $m$, or even continuously differentiable as will be proved in subsection \ref{sec5.3} below.
The solution $u(x;\beta,m)$ is given as the inverse function of the following indefinite integrals:
\begin{align}
\dfrac{1}{\sqrt{2}}\int_{u(x)}^{\beta}\dfrac{dw}{\sqrt{\F_{0}(k)-\F_{0}(w)}}&=\sqrt{\gamma}(l-x)\qquad \text{for}\quad 0<x<l,\label{int1}\\
\dfrac{1}{\sqrt{2}}\int_{\beta}^{u(x)}\dfrac{dw}{\sqrt{\F_{1}(p)-\F_{1}(w)}}&=\sqrt{\gamma}(x-l)\qquad\text{for}\quad l<x<1.\label{int2}
\end{align}


\subsection{Boundary layer and interior transition layer}
In this section we prove assertions (i),
(ii) of Theorem \ref{MAIN}.
First,
we
 consider the behavior of $u(x;\beta,m)$ and $\gamma(\beta,m)$ as $m\downarrow 0$.
 From $\F_{0}(k)=m^{2}/2+\F_{0}(\beta)$,
 we see immediately that $k\to \beta$ as $m\downarrow 0$.
 Similarly,
 $p\to \beta$ as $m\downarrow 0$. 
 This means that $u(x;\beta,m)$ converges to $\beta$ uniformly in $x$ as $m\downarrow 0$.

 Observe that by the continuity of $F_{0}(u)$ and $F_{0}(\beta)<0$,
 there exists $0<k_{0}<\beta$ such that $F_{0}(w)\leq F_{0}(\beta)/2$
 for any $w\in [k_{0},\beta]$.
 Therefore,
 by the mean value theorem,
 we have
$$
\F_{0}(k)-\F_{0}(w)=-F_{0}(k+\theta(w-k))(w-k)\geq -\dfrac{1}{2}F_{0}(\beta)(w-k)
$$
 for $k_{0}\leq k\leq w\leq \beta$.
 Thus,
 \begin{align*}
 M=\int_{k}^{\beta}\dfrac{dw}{\sqrt{2(\F_{0}(k)-\F_{0}(w))}}&\leq \dfrac{1}{\sqrt{|F_{0}(\beta)|}}\int_{k}^{\beta}\dfrac{dw}{\sqrt{w-k}} 
=\dfrac{2\sqrt{\beta-k}}{\sqrt{|F_{0}(\beta)|}}\to 0
 \end{align*}
 as $k \uparrow \beta$.
 In the same way,
 we can prove that $N\to 0$ as $p\downarrow \beta$.
 We therefore conclude that 
 $$
\gamma(\beta,m)=(M+N)^{2}\to 0\quad \text{as}\quad m\downarrow 0.
$$
 
 Next,
 we prove assertion (ii),
 i.e., we consider the behavior of $u(x;\beta,m)$ and $\gamma(\beta,m)$ as $m\uparrow \min\{\sqrt{2|\F_{0}(\beta)|},\sqrt{2|\F_{1}(\beta)|}\}$.
By applying Lemma \ref{thm:asym-1} in Appendix, we can derive the asymptotic behavior of  $M$ and $N$:
\begin{align}\label{assymptotic1}
M=&\dfrac{1}{\sqrt{|\F_{0}'(0)|}}\log\dfrac{1}{k}+O(1)\quad \text{as}\,\,\,k\downarrow 0,\\ \label{assymptotic2}
N=&\dfrac{1}{\sqrt{|\F_{1}'(u_{1})|}}\log\dfrac{1}{u_{1}-p}+O(1)\quad \text{as}\,\,\,p \uparrow u_{1}
\end{align}
which imply
\begin{align*}
M\to \infty\,\,\text{as}\, \,k\downarrow u_{0}(=0),\,\quad N\to \infty\,\, \text{as} \,\,p\uparrow
 u_{1}.
\end{align*}
 Notice that as $k\to 0$ or $p\to u_{1}$,
it holds that $\gamma(\beta,m)\to +\infty$ by virtue of \eqref{assymptotic1} or \eqref{assymptotic2}.
 We discuss the following three cases separately.
 
 (a)\,Assume that $\F_{0}(\beta)>\F_{1}(\beta)$.
In this case, we have $m \uparrow \sqrt{2|\F_{0}(\beta)|}$, 
 and hence $k \to u_{0}$ but $p$ does not approach $u_{1}$ (see \eqref{houteishiki}).
 Therefore,
from \eqref{assymptotic1}-\eqref{assymptotic2} it follows that $M \to +\infty$ and $N$ remains
bounded as $m \uparrow \sqrt{2|\F_{0}(\beta)|}$.
We thus obtain
\begin{align*} 
 l(\beta,m) = M/(N+M)  \to 1 \quad \text{as} \quad m \uparrow 
\sqrt{2|\F_{0}(\beta)|}.
\end{align*}
We can prove that $u(x;\beta,m) \to 0 $ locally uniformly in $[0,1)$ as $ m \uparrow \sqrt{2 |\F_0(\beta)|} $ in the same way as in the case (b) below.

(b)\,Assume that $\F_{0}(\beta^{\ast})=\F_{1}(\beta^{\ast})$ for some $\beta^\ast$.
Then $m \uparrow \sqrt{2|\F_{0}(\beta^{\ast})|}$ and 
$m \uparrow \sqrt{2|\F_{1}(\beta^{\ast})|}$ simultaneously,
so that $k \to u_{0} =0$ and $p \to u_{1}$. 
Consequently,
by virtue of \eqref{assymptotic1} and \eqref{assymptotic2},
$M \to +\infty$ and $N \to + \infty$ at the same time,
and we have 
\begin{align}\label{5-150}
& l(\beta, m) \notag \\
& = 
\left(
\frac{1}{\sqrt{|\F_{0}'(0)|}} \log \frac{1}{k}
+O(1) \right)
\Bigg/ \left(
\frac{1}{\sqrt{|\F_{0}'(0)|}} \log \frac{1}{k}
+
\frac{1}{\sqrt{|\F_{1}'(u_{1})|}} \log \frac{1}{u_{1} -p}
+ 
O(1) \right).
 \end{align}
Here we observe that 
$k$ and $R$ are dependent on each other,
hence $p$ may be regarded as a function of $k$.

From \eqref{houteishiki} we know that 
\begin{align*}
\F_{1}(p) = (\F_{0}(k) - \F_{0}(\beta)) +\F_{1}(\beta),
\end{align*} 
and we now have  $\F_{0}(\beta^{\ast})=\F_{1}(\beta^{\ast})$;
therefore $\F_{1}(p) = \F_{0}(k)$.
On the other hand,
recalling that 
$\F_{0}(0) = \F_{0}'(0) =0$ and 
$\F_{1}(u_{1}) = \F_{1}'(u_{1}) =0$, 
we get 
\begin{equation*}
\F_{0}(k)
 =\frac{1}{2} F'_{0}(\theta_{0} k) k^{2}
\quad \text{and} \quad 
\F_{1}(p)
 =\frac{1}{2} F'_{1}(u_{1} + \theta_{1}(p-u_{1}) ) (p-u_{1})^{2}
\end{equation*}
for some $0< \theta_{0} <1$
and 
$0< \theta_{1} <1$.
Therefore,
\begin{align*}
F_{1}'(u_{1} + \theta_{1}(p-u_{1}))
(p-u_{1})^{2}
=  F'_{0}(\theta_{0} k) k^{2},
\end{align*}
that is, 
\begin{align*}
u_{1} -p
= k 
\sqrt{
\frac{|F_{0}'(\theta_{0} k)|}
{|F_{1}'(u_{1} + \theta_{1}(p-u_{1}))|}
}
= 
 k 
\sqrt{
\frac{|F_{0}'(0) +o(1)|}
{|F_{1}'(u_{1}) +o(1)|}
}.
\end{align*}
Hence
we have,
as $k \downarrow 0$,
\begin{align} \label{5-151}
\log\frac{1}{u_{1} -p}
= \log \frac{1}{k}
-
\frac{1}{2}
\log \frac{|F_{0}'(0) +o(1)|}
{|F_{1}'(u_{1}) +o(1)|}
=
\log \frac{1}{k}
+O(1). 
\end{align}
Consequently, (\ref{5-150}) and (\ref{5-151}) yield
\begin{align*}
l(\beta,m)
& = 
\frac{
\sqrt{F'_{1}(u_{1})}
}
{
\sqrt{F'_{1}(u_{1})}
+
\sqrt{F'_{0}(0)}
}
+o(1)
\quad 
\text{as} \quad 
m \uparrow |\F_{0}(\beta^{\ast})|.
\end{align*}
Let us put 
$l^{\ast} = \sqrt{|F'_{1}(u_{1})|}/ 
( \sqrt{|F'_{1}(u_{1})|}+  \sqrt{|F'_{0}(0)|})$.
Let $\kappa$ be any number satisfying $0 < \kappa < \beta$.
Since $u(x)$ is strictly increasing, there is a unique 
$x_{\kappa} \in (0, l(\beta,m))$ such that 
$u(x_{\kappa},\beta, m) = \kappa$.
Then by \eqref{int1} we have
\begin{align*} 
\int_{\kappa}^{\beta}
\frac{dw}{\sqrt{2(\F_{0}(k) - \F_{0}(w))}}
=
(M+N)(l(\beta,m) -x_{\kappa}),
\end{align*}
and hence
\begin{align*}
0<l(\beta,m) - x_{\kappa}
= 
\frac{1}{M+N}
\int_{\kappa}^{\beta}
\frac{dw}{\sqrt{2(\F_{0}(k) - \F_{0}(w))}}.
\end{align*}
Note that
\begin{align*}
\int_{\kappa}^{\beta}
\frac{dw}{\sqrt{2(\F_{0}(k) - \F_{0}(w))}}
=O(1) \quad \text{as} \quad k \downarrow 0.
\end{align*}
Therefore,
by virtue of $M+N \to +\infty$,
we obtain 
\begin{align*}
|l(\beta,m) -x_{\kappa}|
\to 0 
\quad \text{as} \quad m \uparrow \sqrt{2|\F_{0}(\beta^{\ast})|}.
\end{align*}
Since $l(\beta,m) \to l^{\ast}$,
we conclude that $x_{\kappa} \to l^{\ast}$ as 
$m \uparrow  \sqrt{2|\F_{0}(\beta^{\ast})|}$.
From the inequality $0<u(x;\beta, m) \leq \kappa$
for $0 \le x \le x_{\kappa}$, we therefore
see that $u(x;\beta,m) \to 0$ as $m \uparrow \sqrt{|\F_{0}(\beta^{\ast})|}$
uniformly on each interval $0 \le x \le l^{\ast} - \varepsilon$,
where $\varepsilon>0$ is any small constant.
In the same way,
we can prove that $u(x;\beta,m) \to u_{1}$ as 
$m \uparrow \sqrt{2|\F_{0}(\beta^{\ast})|}$ uniformly on each interval 
$l^{\ast} + \varepsilon \le x \le 1$.

(c)\, Finally,
assume that $\F_{0}(\beta) < \F_{1}(\beta)$.
Then as 
$m \uparrow \sqrt{2|\F_{0}(\beta^{\ast})|}$,
we have $p \to u_{1}$,
but $k$ does not approach $u_{0} =0$.
Hence by the same method as in (a),
we obtain $l(\beta,m) \to 0$,
and $u(x;\beta,m) \to u_{1}$
locally uniformly in the interval $0<x \le 1$. Hence we finish the proof of the theorem.
\qed

\bigskip

\begin{figure}[th!]
 \begin{minipage}{0.5\hsize}
  \begin{center}
  \includegraphics[width=80mm]{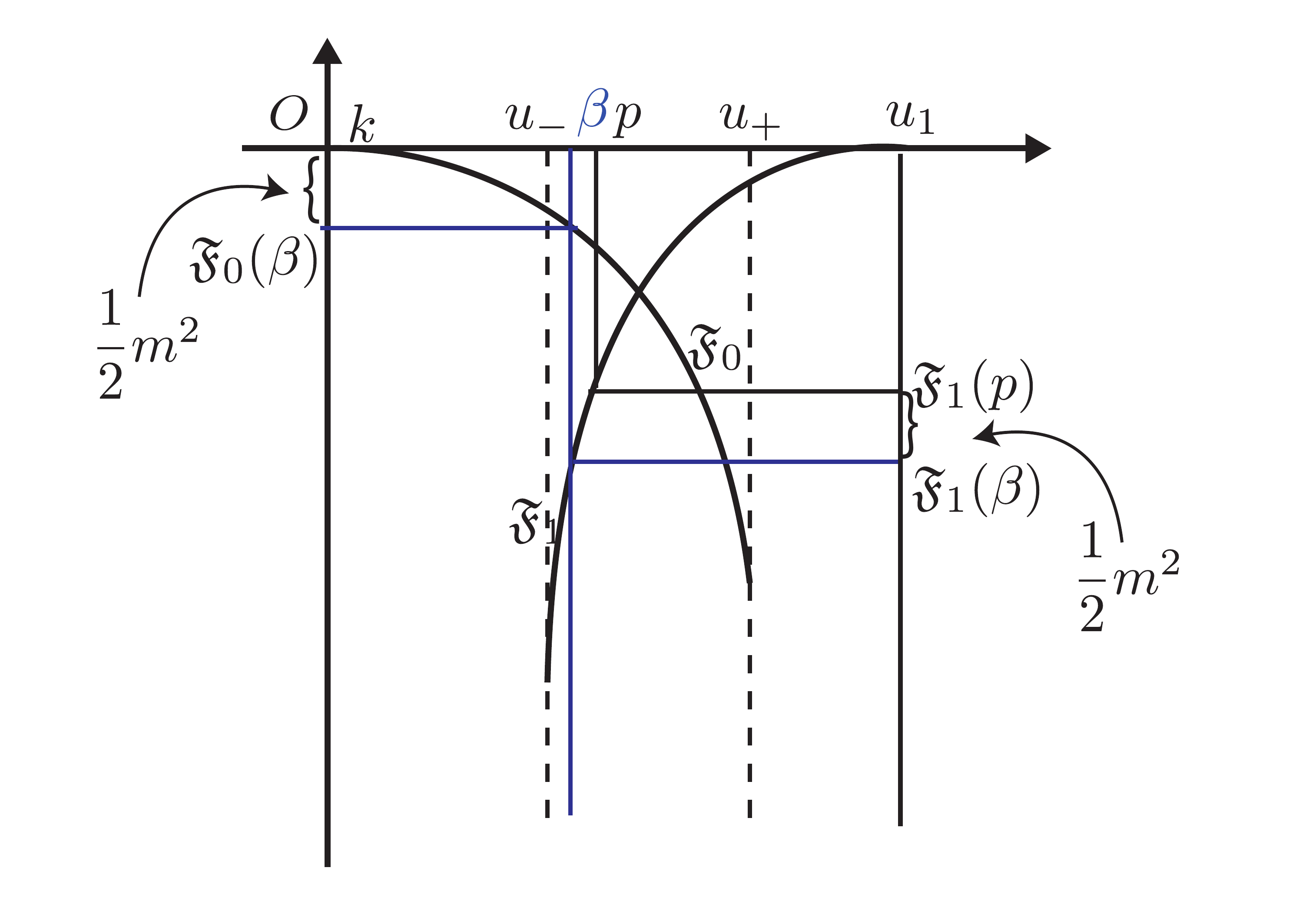}
 \label{a}
  \end{center}
  \vspace {-0.5cm}
  \caption{Case (a)\,\,where ${\mathit\Delta}=m^{2}/ 2$.}
 \end{minipage}
 \begin{minipage}{0.5\hsize}
  \begin{center}
   \includegraphics[width=80mm]{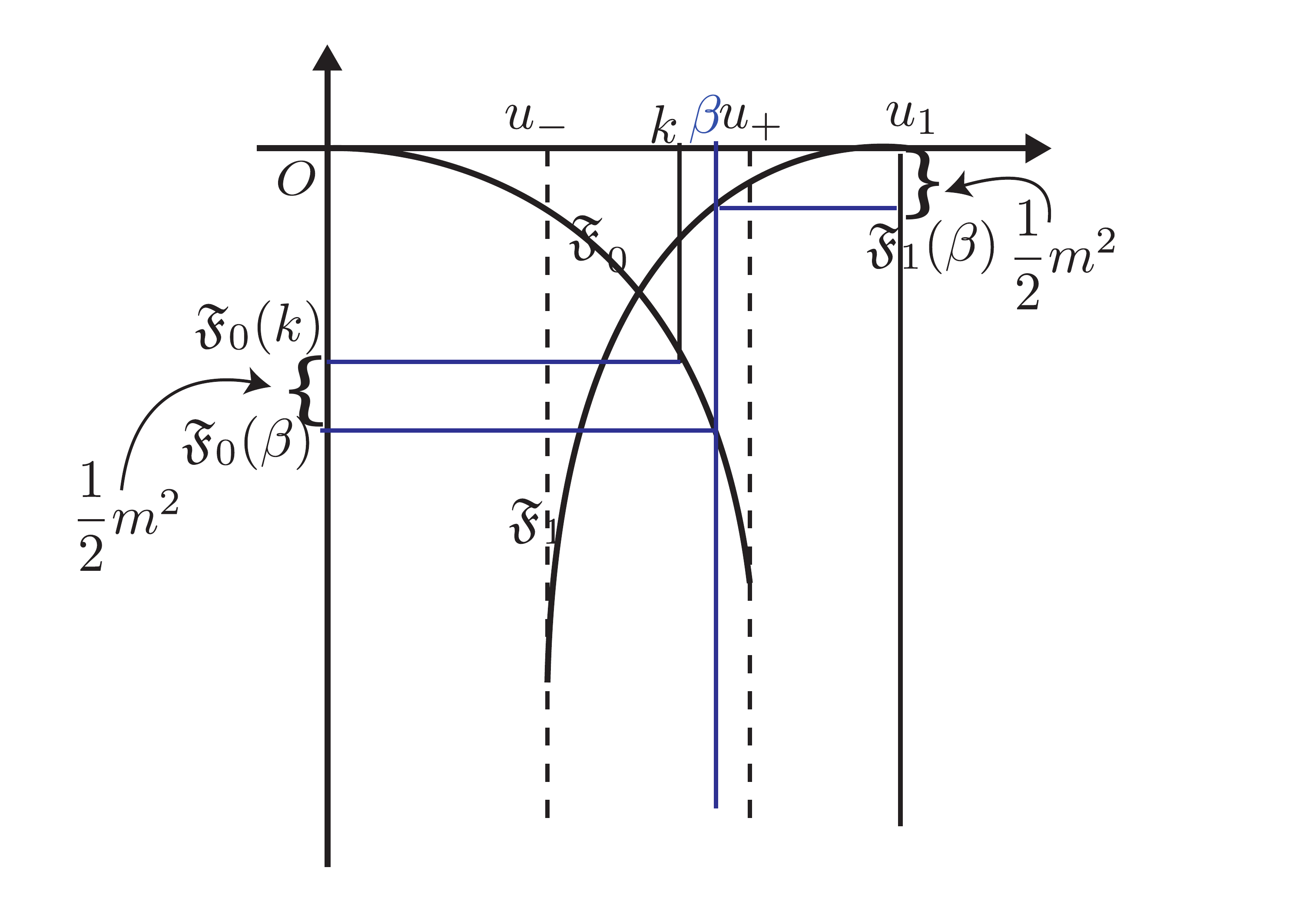}
   \label{c}
  \end{center}
  \vspace{-0.5cm}
  \caption{Case (c)\,\,where ${\mathit\Delta}=m^{2}/2$.}
 \end{minipage}
  \begin{minipage}{0.5\hsize}
   \vspace{1cm}
  \begin{center}
  \includegraphics[width=70mm]{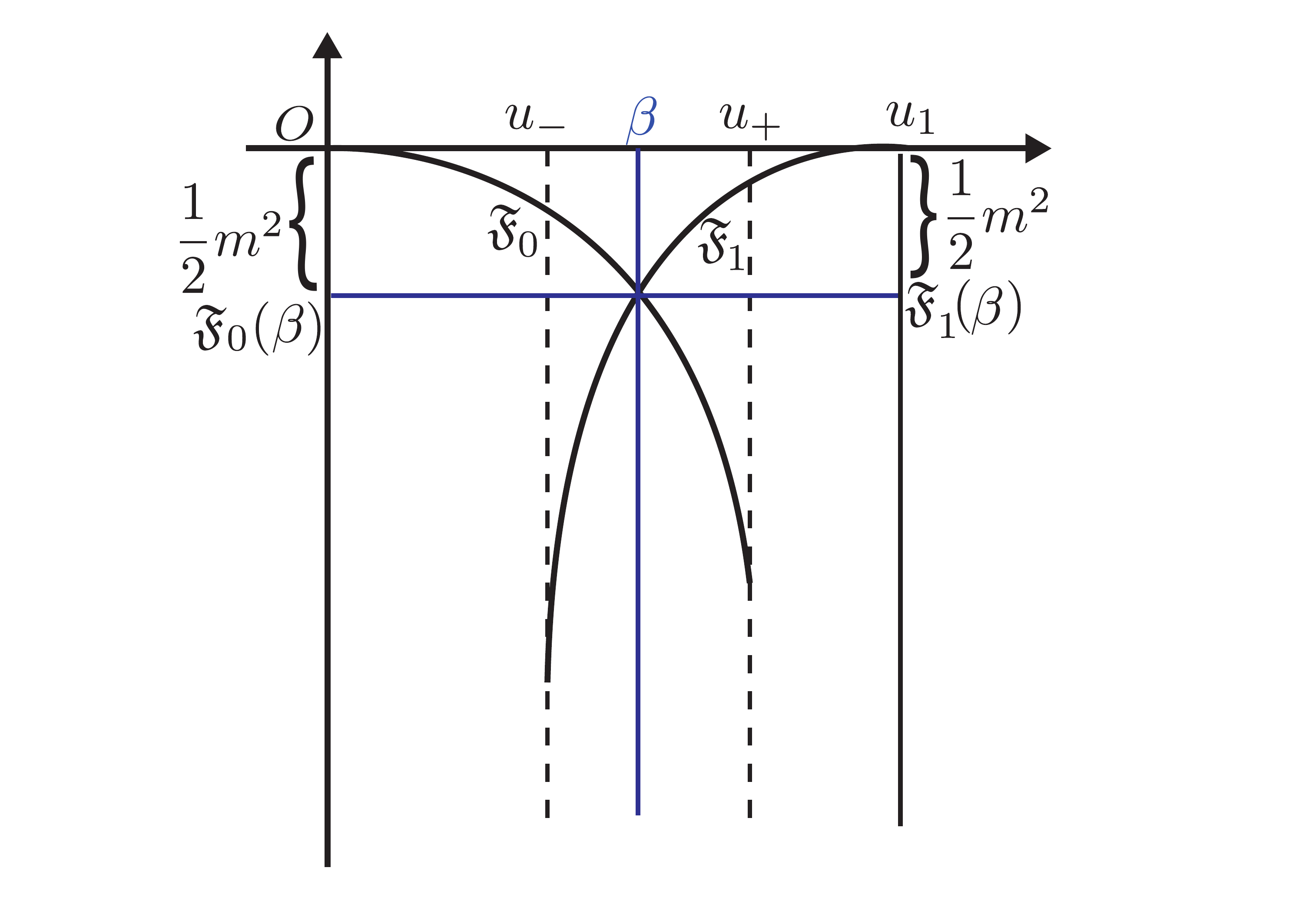}
 \label{b}
  \end{center}
  \vspace{-0.5cm}
  \caption{Case (b)\,\,where ${\mathit\Delta}=m^{2}/2$.}
 \end{minipage}
 \begin{minipage}{0.5\hsize}
  \begin{center}
      \label{}
  \end{center}
 \end{minipage}
\end{figure}


\subsection{A sufficient condition for the uniqueness}\label{sec5.3}
So far we have constructed a one-parameter family of solutions 
$\{ (u(x; \beta,m) , \gamma(\beta,m) )\mid 0<m < \min\{ |\F_{0}(\beta)|,
|\F_{1}(\beta)|
\} \}$ of \eqref{M}.
We know that, for each $\gamma>0$,
there exists at least one $m$ such that $\gamma(\beta,m) = \gamma$. 
In this section,
we give a sufficient condition which guarantees the uniqueness of $m$.
We define functions $M(k)$ and $N(p)$ by 
\begin{align*}
\left\{
\begin{aligned}
& M(k)= 
\int_{k}^{\beta}
\frac{dw}{\sqrt{2(\F_{0}(k) - \F_{0}(w))}}
&&\quad \text{for} \quad 0<k < \beta, \\
& N(p)= 
\int_{\beta}^{p}
\frac{dw}{\sqrt{2(\F_{1}(p) - \F_{1}(w))}}
&&\quad \text{for} \quad \beta<p < u_{1}. 
\end{aligned}
\right.
\end{align*}
Let $p(k)$ be the unique solution of the algebraic equation
\begin{align} \label{1001}
\F_{1}(p) = \F_{0}(k) -\F_{0}(\beta)
+ \F_{1}(\beta).
\end{align}
Then we have 
\begin{align*}
\frac{dp}{dk}
=
\frac{F_{0}(k)}{F_{1}(p(k))}
<0.
\end{align*}
Now we define a function $T(k)$ by 
\begin{align*}
T(k) = M(k) +N(p(k)) \quad \text{for}
\quad 0 < k< \beta.  
\end{align*}
Given $\gamma_{0} >0$,
let $k$ be such that $\gamma_{0} = T(k)^{2}$.
Then,
for this $k$,
we have a unique $p$ satisfying \eqref{1001}.
Once $k$ and $p$ are determined,  $m$ is given uniquely by \eqref{houteishiki};
and hence $u(x)$ for $\gamma_{0}$ is uniquely determined.

Therefore,
we would like to find a sufficient condition for $dT/dk<0$
on the interval $0<k<\beta$.
Observing that 
\begin{align*}
\frac{dT}{dk}
= 
M'(k) + N'(p(k))p'(k)
= 
M'(k)
+ \frac{ N'(p(k))F_{0}(k)}{F_{1}(p(k))},
\end{align*}
we first consider the sign of $M(k)$.
By a change of integration variable $w = k+(\beta- k) t$,
we see that 
\begin{align*}
\sqrt{2} M(k)
=
\int_{0}^{1}
\frac{(\beta-k)dt}
{\sqrt{
\F_{0}(k)
-
\F_{0}(k+(\beta-k) t)
}}.
\end{align*}
From this expression, we obtain, after differentiation of the integrand,
\begin{align*}
\sqrt{2} M'(k)
& 
= -\frac{1}{2(\beta-k)}
\int_{k}^{\beta}
\frac{
2(\F_{0}(k) 
-\F_{0}(w) )
+
(\beta-k)F_{0}(k)
-
(\beta-w)
F_{0}(w)
}
{
(
\sqrt{
\F_{0}(k)
-
\F_{0}(w)
}
)^{3}
} \,
dw.
\end{align*}
Hence,
by introducing the notation 
\begin{align*}
\theta_{0}(u) 
=2\F_{0}(u) - (\beta-u) F_{0}(u),
\end{align*}
we get the formula
\begin{align*}
2 \sqrt{2}\, M'(k)=-\frac{1}{\beta-k}
\int_{k}^{\beta}
\frac{\theta_{0}(k) -\theta_{0}(w)}
{
(
\sqrt{
\F_{0}(k)
-
\F_{0}(w)
}
)^{3}
} \,
dw.
\end{align*}
On the other hand,
we have 
\begin{align*}
\theta_{0}'(u)
& = F_{0}(u) +(\beta-u)F_{0}'(u). 
\end{align*}
Note also that
\begin{align*}
F_{0}'(u)
& = \frac{d}{du}f(u, h_{0}(u)) 
=
f_{u}(u,h_{0}(u))
+ 
f_{v}(u,h_{0}(u))h_{0}'(u) \\
& 
= 
f_{u}(u,h_{0}(u))
+ 
f_{v}(u,h_{0}(u))
\left(
-\frac{g_{u}(u, h_{0}(u))}
{g_{v}(u,h_{0}(u))}
\right). 
\end{align*}
Hence, if we define 
\begin{align*}
\triangle_{0}(u)
=
f_{u}(u,h_{0}(u))
g_{v}(u,h_{0}(u))
- 
f_{v}(u,h_{0}(u)) 
g_{u}(u,h_{0}(u)),
\end{align*}
then
$$
F_0^\prime(u) =\frac{\triangle_{0}(u)} {g_{v}(u,h_{0}(u))}.
$$
Therefore, for $k<u<\beta$,
the estimate 
\begin{align} \label{1002}
\triangle_{0}(u) >0 
\end{align}
implies $ F_0^\prime(u) < 0 $, which in turn yields that
\begin{align*}
\theta_{0}'(u) <0 \quad \text{for} \quad k<u<\beta.
\end{align*}
Thus we conclude that whenever \eqref{1002} holds for $k<u< \beta$,
we have 
\begin{align*}
M'(k) <0
\quad \text{for} \quad k <\beta.
\end{align*}

Second,
we consider the sign of $N(p)$.
By an argument similar to that for $M'(k)$,
we have
\begin{align*}
2 \sqrt{2}\, N(p)
= \int_{0}^{1}
(p-\beta)
\{
\F_{1}(p) -\F_{1}((p-\beta) t + \beta)
\}^{-{1}/{2}}
dt
\end{align*}
and 
\begin{align*}
\sqrt{2} N'(p)
= 
\frac{1}{p-\beta}
\int_{\beta}^{p}
\frac{\theta_{1}(p) - \theta_{1}(w)}
{
( \sqrt{
\F_{1}(p) -\F_{1}(w)
} )^{3} 
}\,
dw,
\end{align*}
where
\begin{align*}
\theta_{1}(u)
=
2 \F_{1}(u) -F_{1}(u)(u -\beta).
\end{align*}
Note that 
\begin{align*}
\theta_{1}'(u)
=
F_{1}(u) -F_{1}'(u)(u- \beta),
\end{align*}
and that 
\begin{align*}
F_{1}'(u)
= 
\frac{\triangle_{1}(u) }
{g_{v}(u,h_{1}(u))}
\end{align*}
where
\begin{align*}
\triangle_{1}(u)
=
f_{u}(u,h_{1}(u))
g_{v}(u,h_{1}(u))
- 
f_{v}(u,h_{1}(u)) 
g_{u}(u,h_{1}(u)).
\end{align*}
Therefore, 
if 
\begin{align} \label{1003}
\triangle_{1}(u)>0 \quad \text{for} \quad \beta < u < u_{+}
\end{align}
then 
\begin{align*} 
F_{1}'(u)<0 \quad \text{for} \quad \beta < u < u_{+},
\end{align*}
and hence 
\begin{align*} 
\theta_{1}'(u)>0 \quad \text{whenever} \quad \beta < u < u_{+}.
\end{align*}
Thus under the condition \eqref{1003},
we see that $N'(p)>0$.
Since $p'(k) <0$,
we conclude that 
\begin{align*}
T'(k) = M'(k) + N'(p)p'(k) <0
\end{align*}
provided that both \eqref{1002} and \eqref{1003} are satisfied.

 If $\triangle_{0} (u)$ changes the sign, 
it is difficult to determine the sign of $M'(k)$.
On the other hand,
if $k$ is sufficiently small,
then we can apply the method used in Lemma 2.5 (p.~225) of \cite{Takagi1986}
to obtain the estimate
\begin{align*} 
M'(k)<0.
\end{align*}
Similarly we have 
\begin{align*}
N'(p) >0
\end{align*}
if $p$ is sufficiently close to $u_{1}$.
Hence,
we get the uniqueness of $k$ such that $T(k)^2 = \gamma$
whenever $\gamma$ is sufficiently large. This proves the second assertion of  Theorem \ref{MAIN}.

%
%

\subsection{Remarks}
By our main results,
we conclude that 
the distribution of $u(x)$ concentrates on the boundary
when the diffusion constant $1/\gamma$ is sufficiently small,
except the special case $\beta=\beta^{\ast}$. 
In the latter case, $u(x)$ develops a sharp transition layer from 
$u = 0$ to $u = u_{1}$ at an interior point $ x= l $ near $ l^{\ast}$ 
when $1/ \gamma$ is sufficiently small.
For other species in the problem (2) we have the following
representation formula for stationary solutions:
\begin{align*}
\begin{aligned}
& u_{1}(x)=\dfrac{m_{1}(\mu_{2}+d)}{\mu_{1}(\mu_{2}+d)+\mu_{2}u(x)}, \qquad 
u_{2}(x)=\dfrac{m_{1}bu(x)}{\mu_{1}(\mu_{2}+d)+\mu_{2}u(x)}, \\
& u_4(x)(=v(x)) =\begin{cases}
h_{0}(u(x))\quad &\text{for}\quad 0<x<l,\\
h_{1}(u(x))\quad &\text{for}\quad l<x<1.
\end{cases}
\end{aligned}
\end{align*}
Therefore, the free receptor $u_{1}(x)$ 
is a monotone decreasing function of $x$, 
while the bound receptor $u_{2}(x)$ is a monotone increasing function of $x$. 
Recall that $h_{j}(u)\,(j=0,1)$ is an increasing function of $u$.
Therefore, the production rate $u_{4}(x)$ of the ligand  is increasing 
in $x$. 

\captionwidth=145mm
\begin{figure}[ht]
 \begin{center}
 \includegraphics[width=120mm]{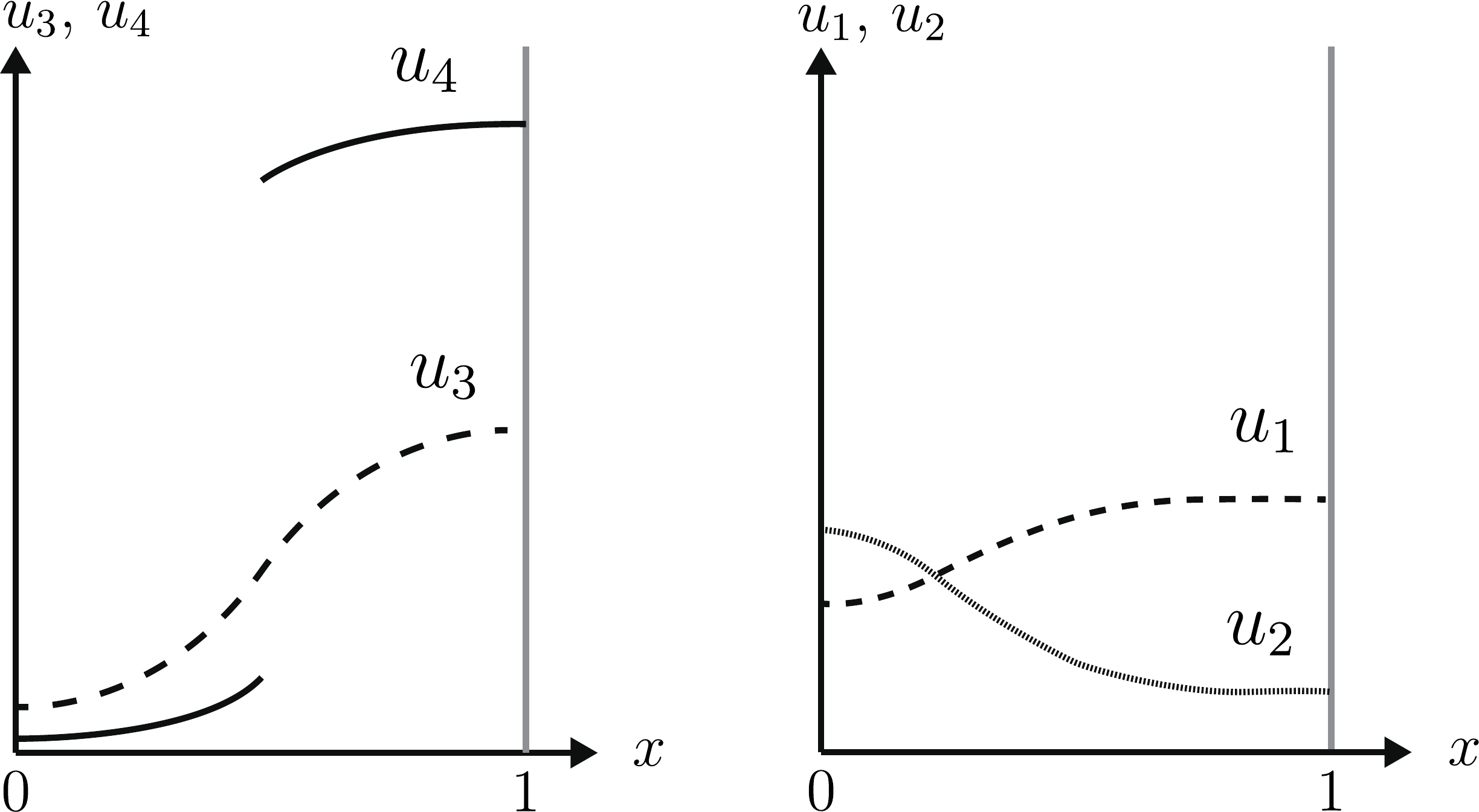} 
 \end{center}
 \hangcaption{
{\footnotesize {\bf Left}:\ Dashed line stands for ligand.\ Solid line stands for production rate of ligand. {\bf Right}:\ Dashed line stands for bound receptor.\ Dotted line stands for free receptor. {\bf Parameters}:\ $D=0.100044,\,\beta=3.731132,\,l=0.46973,\delta=2.5,\,m_{1}=1.5,\,m_{4}=0.75,\,\sigma=0.01,\,\beta_{l}=0.195,\,\mu_{1}=1.0,\,\mu_{2}=1.0,\,\mu_{3}=1.5,\,b=2.0,\,d=1.0$}}
 \label{fig:one}
\end{figure}


We can construct stationary solutions of \eqref{M} which are not monotone increasing.
Let $u_{*}(x)$ be a monotone increasing solution of \eqref{M} with $\gamma=\gamma_{*}$ given by Theorem \ref{MAIN} and let
\begin{align*}v_{*}(x)=
\begin{cases}
h_{0}(u_{*}(x))\quad& \text{for}\quad 0\leq x<l,\\
h_{1}(u_{*}(x))\quad& \text{for}\quad l<x\leq 1.
\end{cases}
\end{align*}
First, we have a monotone decreasing solution $(u_{*}^{-},v_{*}^{-})$ defined by $u_{*}^{-}(x)=u_{*}(1-x),\,v_{*}^{-}(x)=v_{*}(1-x)$.
For each integer $n\geq 2$,
put
$$x_{j}=\dfrac{j}{n}\quad (j=0,1,\dots,n).$$
Define a function $u^{n,+}(x)$ on $0\leq x\leq 1$ by 
\begin{align*}
u^{n,+}(x)=\begin{cases}
u^{+}_{*}(n(x-x_{2k}))\quad &\text{if}\quad x_{2k}\leq x\leq x_{2k+1},\\
u^{-}_{*}(n(x-x_{2k+1}))\quad &\text{if}\quad x_{2k+1}\leq x\leq x_{2k+2}.
\end{cases}
\end{align*} 
and $v^{n,+}(x)$ by
\begin{align*}
v^{n,+}(x)=\begin{cases}
v^{+}_{*}(n(x-x_{2k}))\quad &\text{if}\quad x_{2k}\leq x\leq x_{2k+1},\\
v^{-}_{*}(n(x-x_{2k+1}))\quad &\text{if}\quad x_{2k+1}\leq x\leq x_{2k+2}.
\end{cases}
\end{align*} 
Then $(u^{n,+}(x),v^{n,+}(x))$ is a solution of \eqref{M} with $\gamma=n^{2}\gamma_{*}$.
We call $n$ the {\it mode} of solution.

Similarly,
\begin{align*}
u^{n,-}(x)=\begin{cases}
u^{-}_{*}(n(x-x_{2k}))\quad &\text{if}\quad x_{2k}\leq x\leq x_{2k+1},\\
u^{+}_{*}(n(x-x_{2k+1}))\quad &\text{if}\quad x_{2k+1}\leq x\leq x_{2k+2}.
\end{cases}\\
v^{n,-}(x)=\begin{cases}
v^{-}_{*}(n(x-x_{2k}))\quad &\text{if}\quad x_{2k}\leq x\leq x_{2k+1},\\
v^{+}_{*}(n(x-x_{2k+1}))\quad &\text{if}\quad x_{2k+1}\leq x\leq x_{2k+2}.
\end{cases}
\end{align*} 
gives another solution of mode $n$ of \eqref{M} with $\gamma=n^{2}\gamma_{*}$.

%
%

\appendix

\section{An integral involving a small parameter}

In this appendix we state a lemma on an integral with a small parameter, which was used to compute the principal part of $ M(k) $ and $ N(p) $. The proof may be found in, e.g., Nishiura \cite{nishiura1} for the case $ g \in C^3 $. Here we include an elementary proof that covers the case where $g^{\prime \prime} $ is H\"older continuous.

\begin{lem} \label{thm:asym-1}
Let $g(x)$ be a $C^{2}$ function defined on the interval $0 \le x \le 1$ and satisfies the following assumptions:
\begin{enumerate}[\rm (i)]
\item $g(0)=g'(0)=0$,
\item $g'(x)<0$ \quad \text{and} \quad $g''(x)<0$ for $x\in (0,1)$,
\item 
$g''(x)$ is  H\"older continuous with exponent  $\gamma$ 
 on the interval $0 \le x \le 1$
 $(0 < \gamma <1)$:
\begin{align}
|g''(x)-g''(y)|\leq L|x-y|^{\gamma}\quad \text{for any}\,\,x,y\in[0,1].
\end{align}
\end{enumerate}
Define 
$I(a)$ by 
\begin{align*}
I(a)=\int_{a}^{1}\dfrac{dx}{\sqrt{2(g(a)-g(x))}}\quad \text{for}\,\, a\in (0,1].
\end{align*}
Then
\begin{align}
I(a)=\dfrac{1}{\sqrt{|g''(0)|}}\log{\dfrac{1}{a}}+O(1)\quad \text{as}\,\, a\downarrow
 0.
\end{align}
\end{lem}

We split the integral into two parts and estimate them separately:

\begin{lem} \label{Lem-1}
Let $0<\delta<1$ be fixed and 
set 
\begin{align*}
K(a)=\int_{\delta}^{1}\dfrac{dx}{\sqrt{2(g(a)-g(x))}}\qquad (0< a<\delta<1).
\end{align*}
Then 
$$
K(a)=O(1) \qquad \text{as}\,\,\downarrow 0.
$$ 
\end{lem}

\begin{lem} \label{Lem-2}
For $\delta \in (0, 1) $, let
\begin{align}
J(a)=\int_{a}^{\delta}\dfrac{dx}{\sqrt{2(g(a)-g(x))}}\qquad\text{for}\quad 0<a<\delta.
\end{align}
Then
\begin{align}
J(a)=\dfrac{1}{\sqrt{|g''(0)|}}\log\dfrac{1}{a}+O(1)\qquad\text{as}\quad a\downarrow
 0.
\end{align}
\end{lem}

\noindent{\bf Proof of Lemma \ref{Lem-1}.}
Since $ g^\prime(x) < 0$, we see that $ g(\delta) < g(a) $ if $ a < \delta $, hence $ 2(g(a)-g(x)) > 2(g(\delta) - g(x)) = 2g^\prime(\delta)(x-\delta) - g^{\prime \prime}(\delta + \theta(x-\delta)) (x- \delta)^2 \geq - 2g^\prime(\delta)(x-\delta) $, where $ \theta \in (0,1) $ and we have used the assumption $ g^{\prime \prime}(x) \leq 0 $. Therefore,
$$
K(a)\leq \int_\delta^1 \frac{dx}{ \sqrt{-2g^\prime(\delta) (x-\delta) } } = \frac{ 2(1-\delta) }{ \sqrt{ | g^\prime(\delta) | } },
$$
as desired.
\qed
\bigskip

\noindent {\bf Proof of Lemma \ref{Lem-2}.}
We make a change of integration variable by 
\begin{align} \label{def-xi}
\xi=\sqrt{2(g(a)-g(x))}.
\end{align}
As $x$ increases from $a$ to $\delta$,
the new variable $\xi$ increases from $0$ to $c=\sqrt{2(g(a)-g(\delta))}$. Since the right-hand side of (\ref{def-xi}) is strictly increasing in $x$, for each $ \xi \in [0, c] $, there exists a unique value $ x = x(\xi) $ for which $ \xi = \sqrt{ 2(g(a) - g(x(\xi)) } $ holds. Note also that $ d\xi / dx = - g^\prime(x) / \xi $. Hence,
$$
J(a) = \int_0^c \frac{ d\xi}{- g^\prime(x(\xi)) }.
$$
To proceed further, we solve the equation (\ref{def-xi}) for $x$ in the following way: First, we solve the equation
\begin{align} \label{A6}
-2g(x)  = Ay^2,
\end{align}
where $y$ is a small nonnegative number and $ A=|g^{\prime \prime}(0) |$.
By assumptions (i) and (ii), equation (\ref{A6}) has a unique solution $ x= x_1(y) $ for each $ y \geq 0 $ sufficiently small. Note that $x_1(0)=0 $. Taking $ y = \sqrt{(\xi^2 - 2 g(a) )/A } $ recovers the original equation (\ref{def-xi}).

Since $ g(0)=g^\prime(0)=0 $, we see that $ 2 g(x) = -A(1+\phi(x))x^2 $ with $ \phi(x) = O(x^\gamma) $ as $ x \downarrow 0 $. (In fact, $ \phi(x) = (g^{\prime \prime}(0) - g^{\prime \prime}(\theta x) ) / A $ for some $ \theta \in (0, 1) $.) Then (\ref{def-xi}) is equivalent to the equation $ A(1+\phi(x)) x^2 = Ay^2 $. Since $ x_1(y) $ satisfies this relation, we see that
\begin{align} \label{A7}
(1+\phi(x_1(y)))x_1(y)^2 = y^2.
\end{align}
Form this it follows that $ x_1(y) = (1+o(1))y $ because of $ \phi(x) = o(1) $. We claim that
\begin{align} \label{A8}
x_1(y) = y + O(y^{1+\gamma}) \qquad \text{as} \quad y \downarrow 0.
\end{align}
To see this, we put $ x_1(y) = (1+\eta(y)) y $ and substitute this in (\ref{A7}), obtaining
$$
\{ 1 + \phi((1+\eta(y)) y) \} y^2 (1+\eta(y))^2 = y^2.
$$
This simplifies to 
$$
(2+\eta(y)) \eta(y) = - \frac{ \phi((1+\eta(y))y) } { 1+ \phi((1+\eta(y))y) }.
$$
Recall that $\eta(y) = o(1) $ as $ y \downarrow 0 $. Hence we conclude that $ \eta(y) = O(y^\gamma) $ as $ y \downarrow 0 $.
%
Therefore, $ x(\xi) = x_1( \sqrt{(\xi^2 - 2g(a))/A} ) = \sqrt{ (\xi^2-2g(a))/A } + O\left((\xi^2 - 2g(a))^{(1+\gamma)/2} \right) $.
In view of the fact that $ g^\prime(x) = g^{\prime \prime}(\tilde \theta x) x $ for some $\tilde \theta \in (0,1) $, we have
$$
- g^\prime(x(\xi)) = Ax(\xi) + O( x(\xi)^{1+\gamma}) = \sqrt{ A(\xi^2 - 2g(a)) } \left( 1 + O\left( (\xi^2 - 2g(a))^{\gamma/2} \right) \right).
$$

Consequently,
\begin{align*}
J(a) & = \int_0^c \frac{ d\xi }  { \sqrt{ A(\xi^2 - 2g(a)) } (1+ O( (\xi^2-2g(a) )^{\gamma/2}) ) } \\
& = \int_0^c \frac{ d\xi }{ \sqrt{ A(\xi^2-2g(a)) } } + \int_0^c O\left( (\xi^2-2g(a))^{(\gamma-1)/2} \right)\,d\xi. 
\end{align*}
Let us denote the first integral on the last side by $ J_1(a) $ and the second by $ J_2(a) $. Put
$$
\kappa = \sqrt{-2g(a)} \qquad\text{and}\qquad \xi = \kappa t.
$$
Then $ \kappa = \sqrt{A+o(1)}\, a $ as $ a \downarrow 0$, since $ 2g(a) = g^{\prime \prime}(\theta a) a^2 $. Hence,
$$
\sqrt{A} J_1(a) = \int_0^{c/\kappa} \frac{dt} {\sqrt{ t^2+1 }} = \log \left( c+ \sqrt{c^2+\kappa^2} \right) - \log \kappa = \log \frac{1}{a} + O(1)
$$
as $ a \downarrow 0 $. Finally, we turn to the estimate of $ J_2(a) $:
\begin{align*}
J_2(a) & \leq C \int_0^c \left( \sqrt{ \xi^2 - 2 g(a) } \right)^{\gamma -1} \,d\xi \\
& = \kappa^\gamma \int_0^{c/\kappa} (1+ t^2)^{(\gamma-1)/2}\,dt = \kappa^\gamma \int_0^1 \frac{dt}{ (1+t^2)^{(1-\gamma)/2} } + \kappa^\gamma \int_1^{c/\kappa} \frac{dt}{(1+t^2)^{(1-\gamma)/2}} \\
& \leq \kappa^\gamma + \kappa^\gamma \int_1^{c/\kappa} \frac{dt}{t^{1-\gamma}} \leq \kappa^\gamma + \frac{c^\gamma}{\gamma}.
\end{align*}
This completes the proof of Lemma \ref{Lem-2}.
\qed

\section*{Acknowledgments}
AM-C was supported by European Research Council Starting Grant No 210680 ``Multiscale mathematical modelling of dynamics of structure formation in cell systems'' and Emmy Noether Programme of German Research Council (DFG).  MN and IT were supported in part by JSPS Grant-in-Aid for Scientific Research (A) \#22244010 ``Theory of Differential Equations Applied to Biological Pattern Formation--from Analysis to Synthesis''. Also, IT acknowledges the support of JSPS Grant-in-Aid for Challenging Exploratory Research \#24654037 ``Turing's Diffusion-Driven-Instability Revisited -- from a view point of global structure of solution sets''.

\end{document}